\documentclass[a4paper]{article}

\usepackage{cite}
\usepackage{geometry}
\usepackage{amsmath,amssymb,amsfonts,amsthm}
\usepackage{fontawesome5,hyperref}
\usepackage{graphicx}
\usepackage{textcomp}
\usepackage{dsfont}
\usepackage{enumitem}
\usepackage{nicefrac}
\usepackage{stmaryrd}
\usepackage[dvipsnames]{xcolor}

\newcommand{\I}{\mathbb{I}}

\newcommand{\K}{\mathbb{K}}

\newcommand{\N}{\mathbb{N}}

\renewcommand{\P}{\mathbb{P}}

\newcommand{\R}{\mathbb{R}}

\newcommand{\bOm}{\mathbf{\Omega}}

\newcommand{\bfA}{\mathbf{A}}
\newcommand{\bfB}{\mathbf{B}}

\newcommand{\bfD}{\mathbf{D}}

\newcommand{\bfF}{\mathbf{F}}

\newcommand{\bfL}{\mathbf{L}}

\newcommand{\bfP}{\mathbf{P}}

\newcommand{\bfR}{\mathbf{R}}
\newcommand{\bfS}{\mathbf{S}}

\newcommand{\bfW}{\mathbf{W}}
\newcommand{\bfX}{\mathbf{X}}
\newcommand{\bfY}{\mathbf{Y}}

\newcommand{\cF}{\mathcal{F}}

\newcommand{\cH}{\mathcal{H}}

\newcommand{\cL}{\mathcal{L}}

\newcommand{\cX}{\mathcal{X}}

\newcommand{\vb}{\boldsymbol{b}}

\newcommand{\vf}{\boldsymbol{f}}
\newcommand{\vg}{\boldsymbol{g}}
\newcommand{\vh}{\boldsymbol{h}}

\newcommand{\vp}{\boldsymbol{p}}

\newcommand{\vs}{\boldsymbol{s}}

\newcommand{\vw}{\boldsymbol{w}}
\newcommand{\vx}{\boldsymbol{x}}
\newcommand{\vy}{\boldsymbol{y}}

\newcommand{\vO}{\boldsymbol{0}}

\newcommand{\valpha}{\boldsymbol{\alpha}}

\newcommand{\vgamma}{\boldsymbol{\gamma}}

\newcommand{\vtheta}{\boldsymbol\theta}

\newcommand{\vmu}{\boldsymbol\mu}

\newcommand{\vphi}{\boldsymbol\phi}

\newcommand{\vpsi}{\boldsymbol\psi}
\newcommand{\vomega}{\boldsymbol\omega}

\newcommand{\vpartial}{\boldsymbol\partial}

\newcommand{\mA}{\boldsymbol{A}}

\newcommand{\mK}{\boldsymbol{K}}

\newcommand{\co}{\mathrm{co}}

\newcommand{\inte}{\mathrm{int}}
\newcommand{\x}{\times}

\newcommand{\dist}{\mathrm{dist}}
\newcommand{\st}{\mathrm{s.t.}}

\newcommand{\LV}{\bfL^{\hspace{-0.12em}V}}
\renewcommand{\epsilon}{\varepsilon}

\newtheorem{thm}{Theorem}
\newtheorem{lem}[thm]{Lemma}
\newtheorem{prop}[thm]{Proposition}
\newtheorem{cor}[thm]{Corollary}
\theoremstyle{definition}
\newtheorem{dfn}[thm]{Definition}
\newtheorem{asm}{Assumption}
\newtheorem{pb}[thm]{Problem}
\theoremstyle{remark}
\newtheorem*{rem}{Remark}
\newtheorem*{ex}{Example}

\usepackage{subcaption}
\usepackage{wrapfig}
\usepackage{algorithm,algorithmic}
\usepackage{tikz}
\usepackage{pgfplots}
\usepackage{pgfplotstable,booktabs}
\usepgfplotslibrary{dateplot,fillbetween}
\usepackage{pgf}
\usepgfplotslibrary{colorbrewer}

\title{Robustly Learning Regions of Attraction from Fixed Data}
\author{Matteo Tacchi$^{1}$, Yingzhao Lian$^{2}$, Colin Jones$^{2}$}

\begin{document}

\maketitle

\footnotetext[1]{Univ. Grenoble Alpes, CNRS, Grenoble INP (Institute of Engineering Univ. Grenoble Alpes), GIPSA-lab, 38000 Grenoble, France. e-mail: matteo.tacchi@gipsa-lab.fr}
\footnotetext[2]{Laboratoire d'Automatique, École Polytechnique Fédérale de Lausanne (EPFL), Switzerland. \\e-mail: \{yingzhao.lian, colin.jones\}@epfl.ch }

\begin{abstract}
    While stability analysis is a mainstay for control science, especially computing regions of attraction of equilibrium points, until recently most stability analysis tools always required explicit knowledge of the model or a high-fidelity simulator representing the system at hand. In this work, a new data-driven Lyapunov analysis framework is proposed. Without using the model or its simulator, the proposed approach can learn a piece-wise affine Lyapunov function with a finite and fixed off-line dataset. The learnt Lyapunov function is robust to any dynamics that are consistent with the off-line dataset, and its computation is based on second order cone programming. Along with the development of the proposed scheme, a slight generalization of classical Lyapunov stability criteria is derived, enabling an iterative inference algorithm to augment the region of attraction.
\end{abstract}

\small
\begin{center}
\textbf{Keywords}
\end{center}

Stability of nonlinear systems, Optimization, Uncertain systems, Machine learning, Robust control.

\begin{center}
\textbf{Acknowledgments}
\end{center}

This work was supported by the Swiss National Science Foundation under the “NCCR Automation” grant number 51NF40\_{}180545.

The work of M. Tacchi was also supported by the French company RTE, under the RTE-EPFL partnership n$^\circ$2022--0225.

 \normalsize
 
\tableofcontents

\section{Introduction}\label{sect:intro}

Stability analysis is a major research topic in control science. For instance, given a closed-loop control system and its current configuration (initial condition), deciding whether the trajectories will converge towards a stable equilibrium or not is a crucial problem that control engineers are required to solve in real time on a regular basis. Among various stability criteria, Lyapunov analysis~\cite{lyapunov1992general} plays a key role in this field. In this framework, stability analysis is reformulated into the search for a Lyapunov function (LF), and all initial conditions of trajectories converging towards the stable equilibrium are gathered in the so-called region of attraction (RoA), often approximated from inside with sublevel sets of LFs.

Lyapunov analysis has been widely studied in model-based and data-driven setups, where the user is assumed to have direct access to the model or a high-fidelity simulator. While a system model or simulator may not always be available, it is usually possible and much simpler to measure a finite number of system responses offline. Therefore, stability analysis methods based on fixed measured data become desirable although often challenging. Motivated by this need, this work studies Lyapunov analysis based on a given finite set of measurements of the system response. The proposed approach can learn a piecewise affine (PWA) LF on a compact set without access to a system model or simulator. The contributions of this work are summarized as follows:
\begin{itemize}
    \item We formulate and prove a Lyapunov inference theorem, which generalizes existing Lyapunov stability analysis methods. This generalization is used to expand a prior inner estimate of the region of attraction to a larger set.
    \item We specify our Lyapunov stability criterion for a PWA LF on a compact set, to make it verifiable locally on this compact set rather than at all points.
    \item We make the Lyapunov criterion robust to all models that are consistent with the measured dataset. This criterion is defined on general bounded evaluation function spaces, and has a convex form for Lipschitz functions.
    \item We develop an algorithm to learn a LF, robust to all models that are consistent with the measured dataset. The proposed algorithm only needs to solve a convex second-order cone program, regardless of the underlying unknown dynamics' (non-)linearity. 
    \item We discuss numerical results and properties of the proposed algorithm, mostly focusing on the improvement of its computational efficiency.
\end{itemize}

\subsection*{Previous Work}
Nonlinear Lyapunov stability analysis has been widely studied, where model-based approaches and sampling-based approaches form two main categories. In each approach, a Lyapunov candidate is optimized or synthesized by verifying the Lyapunov stability conditions. In model-based approaches, the knowledge of the underlying model is explicitly used in the search of the LF. In contrast, data-driven approaches train an LF by penalizing the violation of the Lyapunov stability conditions on a dataset. Even though a standard Monte-Carlo sampling scheme can also give a probabilistic guarantee~\cite{hertneck2018learning}, it is always preferable to give a strict qualification in stability analysis. In this case, model-based and most data-driven approaches require explicit knowledge of the model. In particular, verification of the Lyapunov stability condition usually resorts to nonlinear optimization or satisfiability modulo theory (SMT) solvers, such as dReal~\cite{gao2013dreal}. Note that when smooth dynamics are considered,  one can write the Lyapunov stability condition with respect to any Lyapunov candidate into an explicit algebraic form (see e.g \cite{chang2019neural,dai2020counter,kapinski2014simulation}). The SMT solver is accordingly used to check whether these algebraic inequalities are satisfied up to some user-defined tolerance~\cite{gao2012delta}.

To the best of the authors' knowledge, the first numerical method to find an LF solves the Zubov equation~\cite{zubov1964methods}, which models an LF as the solution to a linear partial differential equation (PDE). The approximation of this PDE is then solved via a series expansion~\cite{zubov1964methods}, a collocation method~\cite{giesl2007construction}, etc. One main advantage of the model-based approach is that the a-priori knowledge about the model can be used to reformulate the Lyapunov learning problem into a simpler problem. When polynomial dynamics are considered, a sum of square (SoS) programming relaxation can be used to search for polynomial LFs~\cite{parrilo2000structured}
. Due to the nice algebraic properties of polynomials, the SoS framework has been further used to find the region of attraction~\cite{oustry2019inner,henrion2013convex} and its sparsity structure has been used to improve its scalability~\cite{tacchi2020approximating,ahmadi2019dsos}. Parallel to the studies in polynomial dynamics, PWA dynamics are another active area of research interest~\cite{sun2010stability,lin2009stability}, which comes as a result of the ubiquitous appearance of PWA functions in various controllers, such as ReLU-neural-network-based controllers and linear MPC~\cite{alessio2009survey}. For the PWA setup, optimization based approaches play a central role, such as linear matrix inequalities~\cite{johansson1997computation,henrion2013convex,ravanbakhsh2019learning} and mixed integer programming~\cite{schwan2022stability}. 

Unlike model-based approaches, sampling-based methods rely on an efficient strategy for generating informative samples. Counter-example guided inductive synthesis (CEGIS)~\cite{solar2006combinatorial,solar2008program} is a commonly applied concept in many sample-based approaches (see e.g.~\cite{chen2021learning,ravanbakhsh2019learning,abate2020formal}) when direct access to the model or its simulator is available. During the learning process, they iteratively augment the sample dataset by adding counter examples to the Lyapunov candidate proposed in the current iteration. These algorithms train the LF by penalizing the violation of the Lyapunov stability condition on the samples, and they converge when no further counter example can be generated~\cite{chen2021learning,dai2020counter}.

The search for an LF is usually confined to a specific function class, such as a generalized quadratic form~\cite{johansen2000computation} or a positive definite kernel regressor~\cite{giesl2016approximation}. In this work, we will focus on PWA LFs defined on a compact set. Besides the advantages mentioned in the model-based approach paragraph, PWA Lyapunov candidates have shown nice interplay with Lipschitz dynamics. In particular, when the samples of system dynamics are assigned to the vertices of a grid, a robust Lyapunov stability condition on each simplex can be verified by only considering a tightened Lyapunov condition defined on its vertices. This family of methods is called the continuous piece-wise affine (CPA) method~\cite{marinosson2002lyapunov,julian1999high}. The CPA method has been extended to more general problem setups: differential inclusions~\cite{baier2012linear}, switched systems~\cite{hafstein2007algorithm}, etc. In this work, we also consider Lipschitz dynamics, but we do not assume that the data are located on the vertices. Therefore, we do not term our method a CPA method to avoid unnecessary confusion. A more detailed comparison is presented in Section~\ref{sect:discussion}.

Our method can also be related to~\cite{boffi2021certificates}, which leverages partial knowledge on the system dynamics (such as bounds on its Lipschitz constant) to learn candidate LFs from data and provide probabilistic guarantees for this candidate to generalize to unobserved states. Similarly, we use such knowledge to certify Lyapunov conditions globally, although we choose a robust certification rather than probabilistic, as we aim at avoiding any false positives. Moreover, while~\cite{boffi2021certificates} learns stability certificates, it does not explicitly provide methods for computing the region of attraction, which is the goal of our contribution. Even closer to our work is the reference~\cite{martin2023poly}, where knowledge on the Lipschitz constant is complemented with bounds on higher order derivatives, in order to compute Taylor approximations of the dynamics with deterministic error bounds. This results in polynomial sector conditions and SoS relaxations of Lyapunov constraints. Although less conservative than our method, this framework resorts to semidefinite programming (SDP), while we provide certificates computed from second order cone programming (SOCP) problems.

{The rest of this paper is organized as follows: Section~\ref{sect:pre} states the problem setup. In Section~\ref{sect:theory}, we will first generalize the Lyapunov theorem in Section~\ref{sect:lya}, this generalization will later be used to develop a local Lyapunov condition with PWA Lyapunov candidate in Section~\ref{sect:pwa_lya}. In the sequel, Section~\ref{sect:general_pwa_stb} applies this local condition to a set of uncertain function defined by data, whose robust satisfaction is summarized in Theorem~\ref{thm:robust}. This theorem is later used to define a convex inequality condition for the class of Lipschitz function in Section~\ref{sect:tract}, where the learning problem will be summarized. A comparison between the proposed learning problem and other related works are given in Section~\ref{sect:discussion}. The learnability of the proposed scheme is studied in Section~\ref{sect:valid}, after which the proposed learning problem is recast into an equivalent form to enable higher computational efficiency in Section~\ref{sect:comp_imprv}. The general learning algorithm are summarized in~\ref{sect:algo} with a numerical validation in Section~\ref{sect:result}. A conclusion wraps up this paper in Section~\ref{sect:conclusion}.}

\vspace{1em}

\noindent \textbf{Notation:} $\{\vx_i\}_{i\in\mathcal{I}}$ is a set indexed by $\mathcal{I}$, and when there is no confusion, we drop the index set with $\{\vx_i\}$. $\mathbb{N}_a$ denotes the set of positive integer less than $a$. $\mathbb{R}_+$ denotes the set of non-negative real numbers. $\textbf{0}$ is a zero vector. $\bfB(\vx,r)$ denotes an open euclidean ball centred at $\vx$ with radius $r$. $\co(\bfX)$ denotes the convex hull of $\bfX \subset \mathbb{R}^d$. $\cL_d(\bfX)$ is the set of Lipschitz functions from $\bfX$ to $\R^d$. $\bfX\setminus \bfY:=\{\vx\in \bfX\;|\;\vx\notin \bfY\}$ for all $\bfY \subset \bfX$. $\lVert \vx\rVert := (|x_1|^q+\ldots +|x_d|^q)^{\nicefrac{1}{q}}$ denotes the $q$-norm of $\vx \in \mathbb{R}^d$ for some fixed $q \in [1,\infty]$ (with convention $q=\infty \implies \|\vx\| = \max_{1\leq j \leq d}|x_j|$), whose dual norm is denoted by $\lVert \vx\rVert^*$. $|\bfS|$ denotes the cardinal of a set $\bfS$. If $\bfX$ is a topological space, then $\mathrm{int}(\bfX)$, $\overline{\bfX}$ and $\partial \bfX$ denote its interior, closure and boundary, respectively.

Finally, due to ambiguity in the literature, we indicate our definition of polyhedra and polytopes:
\begin{itemize}
\item A polyhedron $\bfP \subset \mathbb{R}^d$ is an intersection of finitely many half spaces: $\exists c \;\in \mathbb{N}, \mA \in \mathbb{R}^{c\times d}, \vb \in \mathbb{R}^c$ such that
$$ \bfP = \left\{\vx \in \mathbb{R}^d \; | \; \mA\,\vx - \vb \in (\mathbb{R}_+)^c\right\}. $$
\item A polytope is a bounded, finite union of polyhedra (which is not necessarily convex). In particular, a convex polytope is a bounded polyhedron, and is the convex hull of its vertices.
\end{itemize}

\section{Learning regions of attraction from data} \label{sect:pre}

\subsection{Unknown dynamic system, fixed data}\label{sect:stage}
In our problem setup, we consider an \textbf{unknown} continuous time dynamic system of dimension $d$ on a bounded open set $\bfX\subset\R^{d}$:
\begin{subequations} \label{eq:system}
    \begin{align}
    	& \dot{\vx} = \vf(\vx)\; \label{eq:ode} \\
    	& \forall t \geq 0, \quad \vx(t) \in \bfX \label{eq:const}
    \end{align}
\end{subequations}
whose known behavior is summarized as follows:
\begin{asm} \label{asm:regularity} \leavevmode
\begin{enumerate}
\item $\vf$ is Lipschitz continuous, so that the solution to~\eqref{eq:system} with initial condition $\vx(0)=\vx_0$ exists and is unique, and we denote by $\vx(t|\vx_0)$ the corresponding flow at time $t$. 
\item $\vf$ has a single equilibrium point (EP) in $\bfX$, and this EP is in $\vO$ and is locally asymptotically stable (LAS).
\end{enumerate}
\end{asm}

Assumption \ref{asm:regularity} implies that the hypothesis space $\cF$ for the unknown dynamics $\vf$ is at most the space of Lipschitz continuous vector fields: $\cF \subset \cL_d(\bfX)$ with
$$ \cL_d(\bfX) := \left\{\vh: \bfX \to \R^d \; \middle| \; \sup_{\vx\neq\vy} \frac{\|\vh(\vx) - \vh(\vy)\|}{\|\vx-\vy\|} < \infty \right\}. $$
The goal of this work is to algorithmically learn an inner approximation of the RoA for this LAS EP:

\begin{pb}\label{pb:roa} Find the largest possible set of initial conditions of trajectories converging to $\vO$ when $t$ goes to infinity:
$$
	\bfR := \left\{\vx_0 \in \bfX \; \middle| \;  \lim\limits_{t\to\infty} \|\vx(t|\vx_0)\| = 0 \right\}.
$$
\end{pb}

Even when a model is available, such task a is hard in practice, see e.g.~\cite{tacchi2018lyapunov} and the references therein. When there is no model, some alternative information is needed in order to address Problem~\ref{pb:roa}, summarized in Assumptions~\ref{asm:data} and~\ref{asm:uncertainty}.

\begin{asm} \label{asm:data}
A \textbf{fixed} dataset $\bfD$ of size $n$ sampled from system~\eqref{eq:system} is available: $$\bfD := \{(\vx_i,\vf_i = \vf(\vx_i))\}_{i=1}^{n} \subset \bfX \x \R^d.$$
\end{asm}

Assumption~\ref{asm:data} is the only direct observation that we have of the system behavior, and only covers a finite number of configurations, among the infinite possible samples. Hence, we also need some additional knowledge to be able to generalize from this dataset $\bfD$ to the whole set $\bfX \x \vf(\bfX)$. 

\begin{asm} \label{asm:uncertainty}
From the dataset $\bfD$ and side information $(\Sigma)$, a set-valued map $\bfF : \bfX \rightrightarrows \R^d$ exists, such that for all $\vx \in \bfX$, $\bfF(\vx)$ is compact and convex, and it holds that
\begin{equation} \label{eq:uncertainty}
	\vf(\vx) \in \bfF(\vx)
\end{equation}
\end{asm}

\begin{ex} The side information $(\Sigma)$ can take various forms and is essentially a bound on the complexity of the model $\vf$ outside of the datapoints $\vx_i$. Here are some examples of such side information:
\begin{enumerate}
\item Lipschitz bound: $\exists M>0$ such that $\forall \vx,\vy \in \bfX$, $$\|\vf(\vx) - \vf(\vy)\|\leq M \|\vx-\vy\|.$$

In such case, the uncertainty set is trivially given by
$$ \bfF(\vx) = \bigcap_{i=1}^n \bfB(\vf_i, M\|\vx-\vx_i\|). $$
\item Hilbert bounds~\cite{scharnhorst2021robust}: there is a kernel function $\kappa:\bfX\x\bfX \longrightarrow \R$ generating a reproducing kernel Hilbert space (RKHS) $\cH$ of real functions over $\bfX$ such that $\vf = (f_1,\ldots,f_d) \in \cH^d$, as well as upper bounds $M_1,\ldots,M_d > 0$ such that
$$ \forall j \in \N_d, \quad \|f_j\|_\cH \leq M_j. $$
In this setting, it is possible to learn a representer $\widehat{\vf} = (\hat{f}_1,\ldots,\hat{f}_d) \in \cH^d$ for $\vf$ from the dataset $\bfD$, and $\bfF$ can be determined through the following kernel inputs:
\begin{align*}
& \mK := (\kappa(\vx_i,\vx_j))_{i,j=1}^n \in \R^{n\x n}, \\
& \boldsymbol\kappa := \vx\longmapsto (\kappa(\vx,\vx_i))_{i=1}^n \in \R^n, \\
& P := \vx \longmapsto \sqrt{\kappa(\vx,\vx) - \boldsymbol\kappa(\vx)^\top \mK^{-1} \, \boldsymbol\kappa(\vx)},
\end{align*}
for which~\cite{scharnhorst2021robust} provides a new uncertainty set for $\vf$ (denoting $[\pm a] := [-a,a]$ for $a \geq 0$):
$$ \bfF(\vx) := \widehat{\vf}(\vx) + P(\vx)\prod_{j=1}^d \left[\pm\sqrt{M_j^2 - \|\hat{f}_j\|_{\cH}^2}\right]. $$
\item Derivative bounds~\cite{martin2023poly}: $\vf \in C^{k+1}(\bfX)^d$ and for $j \in \N_d$, $\valpha = (\alpha_1,\ldots,\alpha_d) \in \N^d$ such that $|\valpha| := \alpha_1+\ldots+\alpha_d = k+1$,  there is an $M_{j,\valpha} > 0$ verifying
$$ \forall \vx \in \bfX, \quad \left|\frac{\partial^{k+1} f_j}{\partial \vx^{\valpha}} (\vx) \right| \leq M_{j,\valpha},$$
where $\vx^{\valpha} := x_1^{\alpha_1}\cdots x_d^{\alpha_d}$. Here,~\cite{martin2023poly} resorts to degree $k$ Taylor approximations in $\vomega \in \bfX$
\begin{align*}
& \hat{f}_j(\vx) = \sum_{|\valpha| \leq k} c_{j,\valpha} \, (\vx-\vomega)^{\valpha} \\
& \widehat{\vf}(\vx_i) = \vf_i,
\end{align*}
which we summarize into $\widehat{\vf} \in \cF_\bfD^{\vomega} \subset \R[\vx]_k$, to get a polynomial sum-of-squares error bound
$$\|\vf(\vx) - \widehat{\vf}(\vx)\|^2 \leq \sigma_{\vomega}(\vx) := \sum_{j,\valpha} \frac{M^2_{j,\valpha}}{\valpha!} (\vx-\vomega)^{2\valpha},$$
with $\valpha! := \alpha_1!\cdots\alpha_d!$, and end up with the following uncertainty set, parameterized by a finite set $\bOm$ of approximation points :
$$ \bfF(\vx) = \bigcap_{\vomega \in \bOm}\bigcup_{\widehat{\vf} \in \cF_\bfD^{\vomega}} \bfB\left(\widehat{\vf}(\vx), \sqrt{\sigma_{\vomega}(\vx)}\right).$$
It is worth noticing that when $d=1$, the case $k=0$, $\bOm = \{\vx_1,\ldots,\vx_n\}$ recovers our first example as a special case; however, in higher dimension $d$, the two examples are slightly different, the former bounding norms while the latter bounds projections. Also, all these examples can be modified to account for bounded noise $\vf(\vx_i) = \vf_i + \vw_i$, $\P(\vw_i \leq W)=1$.
\end{enumerate}
\end{ex}

\subsection{Some set-theoretical considerations}

We have defined our dynamic system~\eqref{eq:system} only on a bounded set $\bfX$. Indeed, in the three above examples, when $\|\vx\|$ goes to infinity, so does the size of the uncertainty set $\bfF(\vx)$ (because $\|\vx-\vx_i\|$, $P(\vx)$ and $\sigma_{\vomega}(\vx)$ go to infinity), rendering any learning process ineffective. As a result, one can only aim at certifying stability on a bounded subset of the state space. However, in practice, such assumption is without loss of generality, as most control systems are designed with security constraints under the form of bounds on their state variables. In general, outside of these bounds, the integrity of the system at hand is compromised. For these reasons, we only consider bounded admissible state set $\bfX$ in which we assume that the desired RoA $\bfR$ is included: $\bfR \subset \bfX$. We also make an additional assumption on prior knowledge regarding the RoA:
\begin{asm} \label{asm:prior}
A compact subset of the RoA: $\bfA \subset \mathrm{int}(\bfR)$, is already known, with $\vO \in \inte(\bfA)$.
\end{asm}
 
Assumption~\ref{asm:prior} is a slight generalization of Assumption~\ref{asm:regularity} (which yields an $\bfA = \{\vO\}$), required for practical learning of the RoA. Generally speaking, $\bfA$ models a conservative prior about the region of attraction (RoA) of the unknown dynamic system, e.g. deduced from engineering practice. We stress here that such knowledge can be required in a non-trivial form to learn a finitely parameterized LF, either in a model-based or data-driven setting, as discussed in~\cite{ahmadi2011poly, allgower2016exclude}.

\section{Piecewise affine set membership for Lyapunov inference}\label{sect:theory}
In this section, we first try to augment the prior knowledge of attractivity in $\bfA$ to a larger set in Section~\ref{sect:lya}. This result is refined to a Lyapunov candidate from the class of piece-wise affine (PWA) functions in Section~\ref{sect:pwa_lya}. 

\subsection{Lyapunov inference}\label{sect:lya}
Before proceeding to the stability analysis, we introduce three additional concepts on functions $V:\mathbb{R}^{d}\longrightarrow \mathbb{R}$.

\begin{dfn} \leavevmode
\begin{itemize}
\item The strict sub-level set of $V$ with level $a \in \mathbb{R}$ is \[\LV_{a}:=\left\{\vx \in \mathbb{R}^{d} \;|\;V(\vx)< a\right\}.\]
\item The Clarke generalized gradient of $V$ at a point $\vx \in \mathbb{R}^{d}$ is the set given by
$$\vpartial_{\rm C} V(\vx) := \co\left\{
    \vg \in \mathbb{R}^{d} \; \middle| \begin{array}{l} \forall\;\epsilon > 0, \exists \ \vx_\epsilon \in \mathbb{R}^{d} \; \st \\
    \|\vx - \vx_\epsilon\| < \epsilon, \\
    V \text{ is differentiable at } \vx_\epsilon, \\
    \|\vg - \nabla V(\vx_\epsilon)\| < \epsilon
\end{array} \hspace*{-2.35pt} \right\}.$$
\item The Clarke-Lie derivative of $V$ along a set-valued map $\bfF:\R^d \rightrightarrows \R^d$ at a point $\vx \in \R^d$ is
$$ \dot{V}_\bfF(\vx) := \left\{\vf_{\vx}^\top\vg_{\vx} \; \middle| \; \vf_{\vx} \in \bfF(\vx), \; \vg_{\vx} \in \vpartial_{\rm C} V(\vx) \right\}. $$
In particular, if there is an $\vf$ such that for all $\vx$ it holds $\bfF(\vx) = \{\vf(\vx)\}$, we denote
$$ \dot{V}_{\vf}(\vx) := \left\{\vf(\vx)^\top\vg \; \middle| \; \vg \in \vpartial_{\rm C} V(\vx) \right\}. $$
\end{itemize}
\end{dfn}

\begin{rem}
The Clarke gradient is a generalized gradient in the sense that if $V$ is continuously differentiable in a neighbourhood of $\vx$, then trivially $\partial_{\mathrm{C}} V(\vx) = \{\nabla V(\vx)\}$, and if $V$ is convex in a neighborhood of $\vx$, then $\vpartial_{\rm C}V(\vx)$ coincides with the subdifferential of $V$.

In~\cite[Theorem 2.5.1]{clarke} it is proven that if $V$ is Lipschitz continuous in a neighbourhood of $\vx$, then $\partial_{\mathrm{C}} V(\vx) \neq \varnothing$; in such a case, $\vpartial_{\rm C}V(\vx)$ is compact and convex by definition and if $\bfF$ also takes compact and convex values, it trivially follows that $\dot{V}_\bfF$ has the following form:
$$\forall \vx \in \R^d, \; \exists \ell_{\vx} \leq u_{\vx} \; \st \;\dot{V}_\bfF(\vx) = [\ell_{\vx},u_{\vx}].$$
\end{rem}

Then, we introduce a slight generalization of the Krasovsky-LaSalle invariance principle~\cite{la1976stability}.
\begin{lem}\label{lem:stability}
Consider the dynamic system~\eqref{eq:system} and suppose that there is a Lipschitz continuous function $V:\mathbb{R}^{d}\rightarrow \mathbb{R}$ s.t.
\begin{equation} \label{eq:lyapunov}
	\forall \vx \in \bfX \setminus \bfA, \quad 
	\max \dot{V}_{\vf}(\vx) < 0.
\end{equation}

Then, for all $a \in \R$ such that $\LV_a \subset \bfX$, it holds
\begin{equation} \label{eq:convergence}
\forall \vx_0 \in \LV_a, \qquad \lim\limits_{t\to\infty} \dist\left(\vphantom{\sum}\vx(t|\vx_0),\bfA\right) = 0,
\end{equation}
where $\dist(\vp,\bfS) = \inf_{\vs \in \bfS}\|\vp-\vs\|$ denotes the distance between point $\vp \in \R^d$ and set $\bfS \subset \R^d$.
\end{lem}

\begin{rem} The only difference with the original Krasovsky-LaSalle invariance principle is the relaxation of the original continuous differentiability assumption on $V$ into a Lipschitz continuity assumption, and accordingly the use of the Clarke-Lie derivative. This slight modification does not drastically change the proof, thanks to \cite[Lemma 2.15]{MDR2020} which generalizes the chain rule to this setting. As usual, the case $\bfA = \{\vO\}$ covers standard Lyapunov local asymptotic stability theorems (constraints on the values of $V$ being replaced by boundedness of $\bfX$).
\end{rem}

\begin{rem} One can immediately notice that if moreover $\bfA$ is as in Assumption~\ref{asm:prior} (which is not required in Lemma~\ref{lem:stability}), then condition~\eqref{eq:lyapunov} implies that $\vf$ does not vanish outside of the RoA of the LAS EP $\vO$. In other words, $\vO$ is implicitly required to be the \textit{only} equilibrium point in $\bfX$. To relax this implicit condition, a slight modification of Lemma~\ref{lem:stability} can be performed: instead of condition~\eqref{eq:lyapunov}, one can ask for the existence of an $a \in \R$ such that $\LV_a \subset \bfX$ and $\max \dot{V}_{\vf}(\vx) < 0$ only holds for $\vx \in \LV_a \cap (\bfX \setminus \bfA)$ (instead of for all $\vx \in \bfX \setminus \bfA$). Then, conclusion~\eqref{eq:convergence} would hold only for such $a \in \R$. While allowing for the existence of alternative equilibria outside of $\LV_a$, this condition is also more difficult to translate into a convex optimization constraint, hence we did not implement it in our numerical experiments, and focused on the case when $\bfX$ is tailored to contain no EP other than $\vO$.
\end{rem}

In the case of unknown dynamics, condition~\eqref{eq:lyapunov} is impossible to check, as $\dot{V}_{\vf}$ takes unknown values. Hence, in our Lyapunov inference theorem we will replace it with a more conservative but certifiable condition in terms of the uncertainty set $\bfF(\vx)$ defined in~\eqref{eq:uncertainty}.

\begin{thm} \label{thm:stability} Consider the dynamic system~\eqref{eq:system} together with the set-valued map $\bfF$ defined by~\eqref{eq:uncertainty}, and suppose that there exists a Lipschitz continuous function $V:\R^d \longrightarrow \R$ satisfying
\begin{equation} \label{eq:set_lyapunov}
	\forall \vx \in \bfX \setminus \bfA, \quad \max \dot{V}_\bfF(\vx) < 0.
\end{equation}
Then, under Assumptions~\ref{asm:regularity} to~\ref{asm:prior}, for all $a \in \R$ such that $\LV_a \subset \bfX$, it holds $\LV_a \subset \bfR$, i.e.
\begin{equation} \label{eq:stability}
	\forall \vx_0 \in \LV_ a, \qquad \lim_{t\to\infty} \|\vx(t|\vx_0)\| = 0.
\end{equation}
\end{thm}
\begin{proof} It follows from Assumption~\ref{asm:prior} (compact $\bfA$ included in open $\mathrm{int}(\bfR)$) that there exists a $\delta > 0$ such that, if $\dist(\vx,\bfA) < \delta$, then $\vx \in \bfR$\footnote{take e.g. $\delta = \min_{\vx \in \bfA, \vy \in \partial\bfR}\|\vx-\vy\|$; the $\min$ is well-defined as the minimum of a continuous function over the compact set $\bfA \x \partial\bfR$, whose positivity is given by $\bfA \cap \partial \bfR = \varnothing$.}. Existence of a $\tau > 0$ such that $$\dist\left(\vphantom{\sum}\vx(\tau|\vx_0),\bfA\right) < \delta$$ follows from the fact that for all $\vx \in \bfX$, $\vf(\vx) \in \bfF(\vx)$, so that $\dot{V}_{\vf}(\vx) \subset \dot{V}_\bfF(\vx)$ and hence $\max \dot{V}_{\vf}(\vx) \leq \dot{V}_\bfF(\vx) < 0$ recovers~\eqref{eq:lyapunov} from~\eqref{eq:set_lyapunov} and we can apply Lemma~\ref{lem:stability}. As a result, using the semigroup property, it holds
\begin{align*}
    \lim_{t\to\infty} \|\vx(t|\vx_0)\| & = \lim_{t\to\infty} \|\vx(t+\tau|\vx_0)\| \\
    & = \lim_{t\to\infty} \|\vx(t|\vx(\tau|\vx_0))\| = 0.
\end{align*}
\end{proof}

\begin{rem} The point for working with a neighborhood $\bfA$ of the equilibrium in Theorem~\ref{thm:stability} is that this allows us to exclude $\bfA$ from the set on which the Lyapunov condition~\eqref{eq:set_lyapunov} has to be checked. Indeed, on such neighborhood the function $\vf$ becomes vanishingly small, so that in practice, enforcing condition~\eqref{eq:lyapunov} on a finitely parameterized Lyapunov candidate is more challenging (if not impossible) in this region of the state space, see e.g.~\cite{ahmadi2011poly, allgower2016exclude}.

Moreover, similarly to Lemma~\ref{lem:stability}, Theorem~\ref{thm:stability} also includes an implicit condition: in order to allow for the existence of a Lyapunov candidate $V$ satisfying condition~\eqref{eq:set_lyapunov}, i.e. 
$$\forall \vx \in \bfX \setminus \bfA, \qquad 0 > \max\dot{V}_\bfF(\vx) = \max_{\substack{\vf \in \bfF(\vx) \\ \vg \in \vpartial_{\rm C}V(\vx)}} \vf^\top\vg ,$$
$\bfF(\vx)$ should not contain $\vO$ when $\vx \in \bfX \setminus \bfA$, because if $\vO \in \bfF(\vx)$ then $0 \in \dot{V}_\bfF(\vx)$ and condition~\eqref{eq:set_lyapunov} cannot hold, regardless of $V$. In particular, considering the fact that $\bfF$ models uncertainty over $\vf$, it is reduced to a point only on the finite dataset, where $\bfF(\vx_i) = \{\vf_i\}$ (assuming perfect measurements), and grows in size with the distance to datapoints. Hence, it is very likely that on a neighbourhood of $\vO$, and given that $\bfF(\vO) = \{\vO\}$, it will hold $\vO \in \bfF(\vx)$ (this happens for example in the very simple case of a Lipschitz bound $M > 1$). This is another justification for exluding a whole neighbourhood $\bfA$ of $\vO$ in condition~\eqref{eq:set_lyapunov}.
\end{rem}

\subsection{Piecewise affine Lyapunov function}\label{sect:pwa_lya}
PWA functions have strong modelling capability because they are dense in the space of continuous functions with a compact domain~\cite[Chapter 7.4]{royden1988real}. This section will refine Theorem~\ref{thm:stability} to Lipschitz continuous PWA Lyapunov candidates. For the sake of simplicity, we further assume
\begin{asm}\label{asm:poly}
$\overline{\bfX}$ and $\bfA$ are polytopes.
\end{asm} 
When a set is not a polytope, it can be inner-approximated by a polytope up to arbitrary accuracy, thus this assumption will not limit the application of the proposed analysis. Recall that the definition of polytope used in this paper is not necessarily convex, but rather is a finite union of polyhedra  (Section~\ref{sect:stage}).

We now introduce our Lyapunov candidate under the form of a Lipschitz continuous PWA function. We first consider an $m$-piece tessellation of $\bfX \setminus \bfA$:
\begin{subequations} \label{eq:tessellation}
\begin{align}
	& \overline{\bfX \setminus \bfA} = \bigcup_{k=1}^m \bfY_k \label{eq:cup} \\
	& k\neq k' \Longrightarrow \mathrm{int}(\bfY_k)\cap\mathrm{int}(\bfY_{k'}) = \varnothing \label{eq:cap}
\end{align}
\end{subequations} 
where the $\bfY_k$ are convex polytopes (in addition we define $\bfY_0 = \bfA \ni \mathbf{0}$, not necessarily convex). For $k \in \{1,\ldots,m\}$, we denote the vertices of $\bfY_k$ by $\{\vy_{k,l}\}_{l=1}^{\nu_k}$ ($\nu_k$ hence denoting the number of vertices of $\bfY_k$). Using this structure, a PWA Lyapunov candidate $V$ is defined on $\bfX$ by
\begin{subequations}\label{eq:cpwalf}
\begin{align}\label{eqn:pce-aff}
    \forall k \in \{0,\ldots,m\}, \; \vx \in \bfY_k, \quad V(\vx) := \vg_k^\top \vx+b_k\;.
\end{align}
Obviously, our Lyapunov candidate $V$ should be continuous on $\bfX$: for any common vertex $\vy \in \bfY_{k} \cap \bfY_{k'}$ (i.e. $\exists\; l \in \mathbb{N}_{\nu_k},\; l' \in \mathbb{N}_{\nu_{k'}}$ such that $\vy = \vy_{k,l} = \vy_{k',l'}$), the condition
\begin{align}\label{eqn:pce-aff-cont}
    (\vg_k - \vg_{k'})^\top \vy = b_{k'} - b_k
\end{align}
\end{subequations}
should hold. Then, $V$ is Lipschitz continuous on $\bfX$ with Lipschitz constant 
$$ \sup_{\vx\neq\vy} \frac{|V(\vx)-V(\vy)|}{\|\vx-\vy\|} = \max_{0\leq k \leq m} \|\vg_k\| < \infty. $$

\begin{rem} \label{rem}
Regardless of the tessellation~\eqref{eq:tessellation}, system~\eqref{eq:cpwalf} always admits $\{\vg_k,b_k\}_{k=0}^m = \vO$ as a solution. More specific tessellation methods are discussed later in the paper for which non-trivial solutions exist.
\end{rem}

$V$ can then be extended to the whole of $\R^d$ while keeping the same Lipschitz constant, using Kirszbraun's theorem~\cite{kirszbraun}.
Regarding our PWA Lyapunov candidate, the stability condition~\eqref{eq:set_lyapunov} in Theorem~\ref{thm:stability} can be restated as follows:

\begin{lem} \label{lem:pwa-gradient}
Given the tessellation~\eqref{eq:tessellation} of $\bfX$, for $\vx \in \bfX \setminus \bfA$ define the index set 
$$\K(\vx) := \{k \in \N_m \; | \; \vx \in \bfY_k\}.$$
Then, any Lipschitz extension of the function $V$ defined by~\eqref{eq:cpwalf} to the whole $\R^d$ satisfies the following condition:
\begin{equation} \label{eq:pwa-gradient}
	\forall \vx \in \bfX, \quad \vpartial_{\rm C} V(\vx) = \co\{\vg_k \; | \; k \in \K(\vx)\}.
\end{equation}
\end{lem}

\begin{proof}
Let $\vx \in \bfX \setminus \bfA$. If $V$ is differentiable in $\vx$, then for any $k \in \K(\vx)$ it holds
$$ \vpartial_{\rm C}V(\vx) = \{\nabla V(\vx)\} = \{\vg_k\} $$
(which proves that for any $k,k' \in \K(\vx)$, $\vg_k = \vg_{k'}$ and thus by continuity $b_k = b_{k'}$). 
Else, $V$ is not differentiable in $\vx$, and we need the broader definition of the Clarke gradient:
$$ \partial_{\rm C}V(\vx) = \co\left\{ \vg\in\R^d \; \middle| \begin{array}{l}
\forall \epsilon > 0, \exists \vx_\epsilon \in \R^d \; \st \\ \|\vx-\vx_{\epsilon}\| < \epsilon, \\
V \text{ is differentiable in } \vx_\epsilon, \\
\|\vg-\nabla V(\vx_\epsilon)\|<\epsilon
\end{array} \right\}. $$
In our case, there is an $\epsilon > 0$ small enough such that $\|\vx-\vx_\epsilon\| < \epsilon$ implies the existence of a $k \in \K(\vx)$ with $\vx_\epsilon \in \bfY_k$ (because $\bfX$ is open, so that $\vx \notin \partial \bfX$). Then, it follows from the previous argument that if $V$ is differentiable in $\vx_\epsilon$ then $\nabla V(\vx_\epsilon) = \vg_k$. Conversely, for any $k \in \K(\vx)$ there is an $\vx_\epsilon \in \bfY_k \cap \bfB(\vx,\epsilon)$ in which $V$ is differentiable, with again $\nabla V(\vx_\epsilon) = \vg_k$. This allows to recast the Clarke gradient in our special case as
$$ \vpartial_{\rm C}V(\vx) = \co\left\{\vg \in \R^d \; \middle | \begin{array}{l}
\forall \epsilon > 0, \exists k \in \K(\vx) \; \st \\
\|\vg - \vg_k\| < \epsilon
\end{array}\right\}, $$
which is exactly the set announced in~\eqref{eq:pwa-gradient}.
\end{proof}

\begin{prop} \label{prop:stability_pwa}
Let Assumptions~\ref{asm:regularity} to~\ref{asm:poly} hold, and consider the function $V$ defined by~\eqref{eq:tessellation},~\eqref{eq:cpwalf}. If $V$ moreover satisfies the following condition
\begin{equation} \label{eq:pwa-lyapunov}
	\forall k \in \N_m, \, \vx \in \bfY_k, \, \vf \in \bfF(\vx), \quad \vf^\top \vg_k < 0
\end{equation}
Then for all $ a \in \R$ such that $\LV_ a \subset \bfX$, it holds $\LV_ a \subset \bfR$, i.e.
\begin{equation} \tag{\ref*{eq:stability}}
\forall \vx_0 \in \LV_ a, \qquad \lim_{t\to\infty}\|\vx(t|\vx_0)\| = 0.
\end{equation}
\end{prop}

\begin{proof}
Assuming that condition~\eqref{eq:pwa-lyapunov} is satisfied, we want to prove~\eqref{eq:set_lyapunov} in order to apply  Theorem~\ref{thm:stability}. Let $\vx \in \bfX\setminus\bfA$, $\vf \in \bfF(\vx)$, $\vg \in \vpartial_{\rm C}V(\vx)$, and let us prove that $\vf^\top\vg < 0$; closedness of $\dot{V}_\bfF(\vx)$ will conclude. From Lemma~\ref{lem:pwa-gradient}, we have access to a $\vtheta := (\theta_k)_{k \in \K(\vx)} \in [0,1]^{\K(\vx)}$ such that
\begin{align*}
& \sum_{k\in\K(\vx)} \theta_k = 1 \\
& \sum_{k\in\K(\vx)} \theta_k \; \vg_k = \vg
\end{align*}
and hence
$$ \vf^\top\vg = \sum_{k\in\K(\vx)} \theta_k \; \underset{<0}{\underbrace{\vf^\top\vg_k}} < 0. $$
\end{proof}

Finally, we would wrap up this part by sorting out the logic flow in this theorectical Section~\ref{sect:theory} again. The ultimate goal is to extend some prior knowledge of RoA (i.e. $\bfA$) to a larger set $\LV_ a$ via PWA continuous function, which is not smooth. Theorem~\ref{thm:stability} gives this characterization with respect to a Lipschitz continuous Lyapunov candidate via its Clarke gradient evaluation within the set $\bfX\setminus\bfA$. A specific characterization based on a continuous PWA Lyapunov candidate is then summarized in Proposition~\ref{prop:stability_pwa}. This result reformulates the RoA approximation into a negativity test on each cell $\bfY_k$.

\begin{rem}\label{rmk:B_diss}

It is noteworthy that, with a fixed tessellation, the parameters of the Lyapunov candidate on each affine piece (i.e. $\vg_k,b_k$ on $\bfY_k$) can be uniquely determined by the function evaluation on the vertices $\{V(\vy_{k,l})\}_{l=1}^{\nu_k}$.

Another main benefit of a fixed tessellation is that it allows a direct control over the model complexity of the Lyapunov candidate. In particular, consider two Lyapunov candidates $V_1(\vx)$ and $V_2(\vx)$ with their corresponding partitions $\{\bfY_{1,k}\}$ and $\{\bfY_{2,k}\}$. Then, we can state that $V_1(\vx)$ is a refinement of $V_2(\vx)$ (\textit{i.e.} $V_1(\vx)$ has a higher degree of modelling capability than $V_2(\vx)$) if $\forall\;\bfY_{2,k},\;\exists \{\bfY_{1,j}\}_{j\in\mathcal{I}_k}$ such that $\cup_j\bfY_{1,j} = \bfY_{2,k}$. As condition~\eqref{eq:pwa-lyapunov} is local to each cell, if one cell $\bfY_k$ violates the assumptions of Proposition~\ref{prop:stability_pwa}, then we can refine the model locally by further partitioning $\bfY_k$.
\end{rem}

\section{Learning robust Lyapunov functions} \label{sect:lrpwalf}

The developments of the previous section are similar in spirit to set membership methods for dynamical systems, such as~\cite{scharnhorst2021robust,martin2023poly}, although their specific application to computing a region of attraction is relatively new, and their combination with search spaces made of PWA Lyapunov candidates as proposed in Section~\ref{sect:theory} is original. In this section, we further specialize our approach by leveraging specific structures of the uncertainty set, in order to reach unprecedented scalability properties by writing the overall learning process as a second order cone program (SOCP).

Subsection~\ref{sect:general_pwa_stb} studies the negativity condition~\eqref{eq:pwa-lyapunov} for a general hypothesis space, and Subsection~\ref{sect:tract} studies the condition in the most basic hypothesis space, i.e. the Lipschitz function space.

\subsection{Robustly distributing information}\label{sect:general_pwa_stb}

The key to certifying that a PWA Lyapunov candidate is indeed an LF in our setting, consists in verifying the negativity condition~\eqref{eq:pwa-lyapunov}. This condition can be strengthened for $k \in \N_m$ into
$$ v^\star_k := \sup_{\vx\in\bfY_k} \max_{\vf \in \bfF(\vx)} \vf^\top\vg_k < 0, $$
whose left hand side is an optimization problem with linear cost $\vf \longmapsto \vf^\top\vg_k$. Hence, a natural approach would be to consider the feasible set
\begin{equation} \label{eq:feasible}
\bfF(\bfY_k) := \bigcup_{\vx \in \bfY_k} \bfF(\vx)
\end{equation}
and ask two questions about it:
\begin{enumerate}
\item Is $\bfF(\bfY_k)$ compact and convex? If so, then $v^\star_k$ is attained on the set $\bfF(\bfY_k)^\star$ made of its extreme points. 
\item Can $\bfF(\bfY_k)^\star$ be efficiently parameterized? If so, then it may be possible to compute $v_k^\star$ and test for its sign.
\end{enumerate}
If $\bfF:\bfX \rightrightarrows \R^d$ is determined through Hilbert bounds as in~\cite{scharnhorst2021robust}, then convexity of $\bfF(\bfY_k)$ is very difficult to decide. In contrast, if $\bfF$ comes from Lipschitz bounds, then it is possible to prove that $\bfF(\bfY_k)$ is compact and convex. 
However, even in this simple case, the extreme points $\bfF(\bfY_k)^\star$ are difficult to parameterize as soon as $n>1$, unless one notices a specific structure on $\bfF$, which we will now assume.

\begin{asm} \label{asm:localize}
There exists a finite set $\bOm \subset \bfX$ and a family of set-valued maps $\bfF_{\vomega}: \bfX \rightrightarrows \R^d$ indexed by $\vomega \in \bOm$, such that for all $\vx \in \bfX$, it holds
\begin{subequations} \label{eq:localize}
\begin{equation} \label{eq:intersect}
	\bfF(\vx) = \bigcap_{\vomega \in \bOm} \bfF_{\vomega}(\vx).
\end{equation}
and for $\vomega \in \bOm$ the following ``local'' feasible set is compact:
\begin{equation} \label{eq:local_feasible}
\bfF_{\vomega}(\bfY_k) := \bigcup_{\vx \in \bfY_k} \bfF_{\vomega}(\vx).
\end{equation}
\end{subequations}
\end{asm}

\begin{rem} Assumption~\ref{asm:localize} rules out Hilbert bounds from our framework (in their current form). However, it still covers derivative bounds as in~\cite{martin2023poly} as well as Lipschitz bounds.
\end{rem}

An intuitive interpretation of Assumption~\ref{asm:localize} would be that the uncertainty on $\vf(\vx)$ can be sequentially reduced using information related to parameterizing points $\vomega \in \bOm$. From this observation, the information relative to each point $\vomega$ can be distributed in the learning process, using the following result.

\begin{thm} \label{thm:robust} Under Assumption~\ref{asm:localize}, for any $k \in \N_m$,~\eqref{eq:pwa-lyapunov} holds if there exists a map $\vgamma_k : \bOm \longrightarrow \R^d$ such that
\begin{subequations} \label{eq:robust}
\begin{align}
& \sum_{\vomega \in \bOm} \vgamma_{k,\vomega} = \vg_k \label{eq:decomposition} \\
& \sum_{\vomega \in \bOm} \max_{\vf \in \bfF_{\vomega}(\bfY_k)} \vf^\top \vgamma_{k,\vomega} < 0 \label{eq:distribution}
\end{align}
\end{subequations}
\end{thm}

\begin{proof} 
Let $\vomega \in \bOm$, $\vx \in \bfY_k$. Then, by Assumption~\ref{asm:localize}, it holds $\bfF(\vx) \subset \bfF_{\vomega}(\vx)$, so that taking the union over all possible $\vx$ yields $\bfF(\bfY_k) \subset \bfF_{\vomega}(\bfY_k)$. Hence, for any $\vf \in \bfF(\bfY_k)$ it holds
$$ \vf^\top \vgamma_{k,\vomega} \leq \max_{\vf_{\vomega} \in \bfF_{\vomega}(\bfY_k)} \vf_{\vomega}^\top \vgamma_{k,\vomega}, $$
which can be summed over $\vomega \in \bOm$ to get
\begin{eqnarray*}
\vf^\top \vg_k & \stackrel{\eqref{eq:decomposition}}{=} & \sum_{\vomega \in \bOm} \vf^\top \vgamma_{k,\vomega} \\
& \leq & \sum_{\vomega \in \bOm} \max_{\vf_{\vomega} \in \bfF_{\vomega}(\bfY_k)} \vf_{\vomega}^\top\vgamma_{k,\vomega}
\end{eqnarray*}
and we conclude the proof using~\eqref{eq:distribution}.
\end{proof}

\begin{rem}
The key concept behind Theorem~\ref{thm:robust} is the decomposition of the uncertainty set $\bfF(\vx)$ in~\eqref{eq:intersect}. In particular, if $\bOm = \{\vx_i\}_{i=1}^n$ is made of the sampled data points, then the quantity $\max_{\vf_{\vomega} \in \bfF_{\vomega}(\bfY_k)} \vf_{\vomega}^\top \vgamma_{k,\vomega}$ is related to the uncertainty quantified from one data point, which usually has an easy-to-evaluate explicit closed solution. In comparison, the explicit solution is usually not available or difficult to evaluate when the whole dataset $\bfD$ is considered. For example, when considering Lipschitz bounds, the uncertainty boundary quantified by one data point defines a shifted cone. However, the uncertainty upper and lower bounds are PWA and non-trivial to evaluate~\cite{calliess2020lazily} {when the whole dataset $\bfD$ is used}. Moreover, it is also reasonable to consider an RKHS, which underpins various uncertainty quantification methods such as Gaussian process regression~\cite{williams2006gaussian} and deterministic error bound methods~\cite{scharnhorst2021robust}. All these methods require computing the inverse of the Gram matrix or solving a second order cone program, which has an easy-to-evaluate explicit solution only when one data point is considered. The consequence of Theorem~\ref{thm:robust} is that instead of studying compactness and convexity of the $\bfF(\bfY_k)$, one is reduced to studying convexity of the $\bfF_{\vomega}(\bfY_k)$ as well as their extreme points, which can be much simpler.
\end{rem}

\subsection{A Convex tractable case: Lipschitz bounds}\label{sect:tract}

Theorem~\ref{thm:robust} gives a representation of condition~\eqref{eq:pwa-lyapunov}, but such a representation remains abstract and hard to check numerically; for instance, with derivative bounds it holds
$$ \bfF_{\vomega}(\vx) = \bigcup_{\widehat{\vf} \in \cF_\bfD^{\vomega}} \bfB\left(\widehat{\vf}(\vx), \sqrt{\sigma_{\vomega}(\vx)}\right), $$
and although convexity of $\bfF_{\vomega}(\bfY_k)$ can be deduced from the convexity of $\cF_\bfD^{\vomega}$ and the cone of SoS polynomials, its extreme points remain very difficult to identify, so that using Theorem~\ref{thm:robust} does not actually simplify the computations: one still has to solve SDP problems. Thus, we will now recast this representation under a tractable form, in a specific case. Although we have discussed that it is possible to consider a more complex side information on top of Lipschitz bounds in Section~\ref{sect:stage}, this section will show that a convex learning problem exists when we consider this most basic hypothesis. Indeed, a minimal assumption to deduce the general behavior of the dynamics $\vf$ on $\bfX$ from local data $\bfD$, is formulated in terms of bounding the variations of $\vf$.

\begin{asm}\label{asm:lipschitz}
An upper bound $M$ on the Lipschitz constant of the unknown vector field $\vf$ is known:
\begin{equation} \label{eq:lipschitz}
\exists M > 0 ; \quad \forall \vx,\vy \in \bfX, \quad \|\vf(\vx) - \vf(\vy)\|\leq M \|\vx-\vy\|.
\end{equation}
\end{asm}

In this case, the following corollary to Proposition~\ref{prop:stability_pwa} and Theorem~\ref{thm:robust} allows one to drastically simplify the numerical treatment of constraint~\eqref{eq:pwa-lyapunov}.

\begin{thm}\label{thm:rb_stb_Lip}
Let Assumptions~\ref{asm:regularity},~\ref{asm:data},~\ref{asm:prior},~\ref{asm:poly} and~\ref{asm:lipschitz} hold, and consider the function $V$ defined by~\eqref{eq:tessellation},~\eqref{eq:cpwalf}. For $k \in \N_m$, let $\{\vy_{k,l}\}_{l=1}^{\nu_k}$ be the set of vertices of $\bfY_k$, so that 
\begin{subequations} \label{eq:rb_stb_Lip}
\begin{equation} \label{eq:vertices}
\bfY_k = \co\{\vy_{k,l}\}_{l=1}^{\nu_k}.
\end{equation}
Suppose that $\forall k \in \N_m$,~there is a set $\{\vgamma_{i,k}\}_{i=1}^n \subset \R^d$ with
\begin{equation} \label{eq:deco-Lip}
\sum_{i=1}^n \vgamma_{i,k} = \vg_k
\end{equation}
and for all $l \in \N_{\nu_k}$
\begin{equation}\label{eq:robust-lyapunov}
\sum\limits_{i=1}^{n} \vf_i^\top\vgamma_{i,k}+ M \|\vgamma_{i,k}\|^*\| \vx_i - \vy_{k,l}\| <0,
\end{equation}
\end{subequations}
where we recall that there is a $q \in [1,\infty]$, $q':=\nicefrac{q}{(q-1)}$ (with conventions $\nicefrac{1}{\infty}=0$ and $\nicefrac{1}{0} = \infty$) such that for $\vphi \in \R^d$, $\|\vphi\|^q = |\phi_1|^q + \ldots + |\phi_d|^q$ and $\|\vphi\|^{*q'} = |\phi_1|^{q'}+\ldots+|\phi_d|^{q'}$.

Then, $\forall a \in \R$ such that $\LV_a \subset \bfX$, it holds $\LV_a \subset \bfR$, i.e.
\begin{equation} \tag{\ref*{eq:stability}}
\forall \vx_0 \in \LV_a, \qquad \lim_{t\to\infty}\|\vx(t|\vx_0)\| = 0.
\end{equation}
\end{thm}
\begin{proof}
See Appendix~\ref{app:proof}.
\end{proof}

\begin{rem} Theorem~\ref{thm:rb_stb_Lip} exploits two features specific to Assumption~\ref{asm:lipschitz}, namely the fact that $\bfF_i(\vx)$ is centered at $\vf_i$ (in contrast, with Hilbert or derivative bounds, it is centered at some approximant $\widehat{\vf}(\vx)$ of $\vf(\vx)$) and that for all $\vx \in \bfX$, $\max_{\vf \in \bfF_i(\vx)}\|\vf - \vf_i\|$ is upper bounded by a convex function of $\vx$, which attains a maximum on $\bfY_k$ at a vertex $\vy_{k,l}$. In comparison, the shapes of uncertainty sets related to Hilbert or derivative bounds are much more difficult to analyse.
\end{rem}

Theorem~\ref{thm:rb_stb_Lip} implies that, under Assumptions~\ref{asm:regularity},~\ref{asm:data},~\ref{asm:prior},~\ref{asm:poly} and~\ref{asm:lipschitz}, existence of an LF giving access to a robust inner approximation $\LV_a$ of the RoA $\bfR$ of system~\eqref{eq:system}, is obtained as the consequence of the feasibility of the system of constraints~\eqref{eq:tessellation},~\eqref{eq:cpwalf},~\eqref{eq:rb_stb_Lip}.  Further assuming that the tessellation~\eqref{eq:tessellation} is fixed and recalling again that for $\vphi \in \R^d$,
$$ \|\vphi\| = \left(\sum_{j=1}^d |x_j|^q \right)^{\nicefrac{1}{q}}, $$
this corresponds to:
\begin{itemize}
\item linear (polytopic) constraints if $q=1$: $\forall j,k,l$ 
$$ \sum_{i=1}^n \vf_i^\top \vgamma_{i,k} + M \gamma_{i,j,k} \sum_{r=1}^d |x_{i,r} - y_{k,l,r}| < 0 $$
\item second order cone constraints if $q=2$: $\forall k,l$ 
$$ \sum_{i=1}^n \vf_i^\top \vgamma_{i,k} + M \|\vgamma_{i,k}\|\cdot\|\vx_i - \vy_{k,l}\| < 0. $$
\end{itemize}
In contrast, contributions like~\cite{martin2023poly} or~\cite{boffi2021certificates} resort to the more expensive framework of sum-of-squares programming, which belongs to the class of semidefinite programming and whose complexity grows quickly with the state space dimension and degree of the considered certificates, while~\cite{chen2021learning} is based on the nonconvex framework of mixed-integer programming.

Introducing a fixed negativity tolerance $\epsilon > 0$ as well as slack variables $s_{k,l} \in \R$, $k \in \N_m$, $l \in \N_{\nu_k}$, we replace the strict inequality constraint~\eqref{eq:robust-lyapunov} with a large inequality constraint, and learn an LF through solving the following optimization problem.
\begin{pb} \label{pb:lasalle} Find $\{b_k\}_{k=0}^m \subset \R$, $\{\vgamma_{i,k}\}_{1\leq i \leq n}^{1 \leq k \leq m} \subset \R^d$, $\vg_0 \in \R^d$, $\{\{s_{k,l}\}_{l=1}^{\nu_k}\}_{k=1}^m \subset \R$ solution to the program
\begin{subequations} \label{opt:lasalle}
\begin{align}
& s_\epsilon^\star := \min \sum\limits_{k=1}^m \sum\limits_{l=1}^{\nu_k} s_{k,l} \notag \\
& \text{if } \vy_{k,l} \in \bfY_0, \quad \sum_{i=1}^n \vgamma_{i,k}^\top \vy_{k,l} = \vg_0^\top\vy_{k,l} + b_0 - b_k \tag{\ref*{eqn:pce-aff-cont}} \\
& \text{if } \vy_{k,l} \in \bfY_{k'}, \quad \sum_{i=1}^n (\vgamma_{i,k}-\vgamma_{i,k'})^\top \vy_{k,l} = b_{k'}-b_k \tag{\ref*{eq:decomposition}} \\
& \forall k \in \N_m, l \in \N_{\nu_k}, \quad s_{k,l} \geq -\epsilon \label{con:slack-lasalle} \\
& \text{and} \quad \sum_{i=1}^n \vf_i^\top\vgamma_{i,k} + M\|\vgamma_{i,k}\|^*\|\vx_i - \vy_{k,l}\| \leq s_{k,l}. \label{con:lyapunov}
\end{align}
\end{subequations}
\end{pb}

Optimization problem \eqref{opt:lasalle} comes with the following result:

\begin{prop} \label{prop:lasalle}
Under Assumption~\ref{asm:poly} and tessellation given by~\eqref{eq:tessellation}, for any $\vf$ satisfying Assumptions~\ref{asm:regularity},~\ref{asm:data},~\ref{asm:prior}, and~\ref{asm:lipschitz} and any optimal solution to Problem~\ref{pb:lasalle} satisfying
\begin{equation} \label{eq:lasalle}
s^\star_\epsilon = -\epsilon\sum_{k=1}^{m} \nu_k,
\end{equation}
the function $V$ defined for $k \in \N_m$, $\vx \in \bfY_k$ by
\begin{equation} \tag{\ref*{eqn:pce-aff}+\ref*{eq:deco-Lip}}
V(\vx) = \sum_{i=1}^n \vgamma_{i,k}^\top\vx + b_k
\end{equation}
defines a generalized Lyapunov function, such that for all $a \in \R$ satisfying $\LV_a \subset \bfX$, it holds $\LV_a \subset \bfR$, i.e.
\begin{equation} \tag{\ref*{eq:stability}}
\forall \vx_0 \in \LV_a, \qquad \lim_{t\to\infty}\|\vx(t|\vx_0)\| = 0.
\end{equation}
\end{prop}
\begin{proof}
$s^\star_\epsilon = -\epsilon\sum_{k=1}^m \nu_k$ in Problem~\ref{pb:lasalle} directly implies that for all $k \in \N_m$, $l \in \N_{\nu_k}$, it holds $s_{k,l} = -\epsilon$, i.e. constraint~\eqref{con:slack-lasalle} is saturated. Reinjecting in constraint~\eqref{con:lyapunov} yields~\eqref{eq:robust-lyapunov} as $\epsilon > 0$. The conclusion then follows directly from Theorem~\ref{thm:rb_stb_Lip}.
\end{proof}

The implication of this result is strong as it states that under the assumption of Lipschitz dynamics, we can learn/validate a robust estimate of the RoA $\bfR$ through convex programming~\eqref{opt:lasalle} even when the unknown underlying dynamic system $\vf$ is nonlinear.

\begin{rem}
It is possible to adapt the proposed scheme~\eqref{opt:lasalle} to the case where the measurements are contaminated by bounded measurement noise. More specifically, if one has only access to a dataset $\widetilde{\bfD} := \{(\widetilde{\vx}_i, \widetilde{\vf}_i)\}_{i=1}^n \subset (\R^d)^2$ and bounds $W, W'>0$ such that there exists $\bfD = \{(\vx_i,\vf_i)\}_{i=1}^n \subset \bfX\x\R^d$ and a set of random variables $\bfW = \{(\vw_i, \vw_i')\}_{i=1}^n \subset (\R^d)^2$ satisfying $\vf(\vx_i) = \vf_i$ and

$$ \begin{array}{ll}
    \widetilde{\vx}_i = \vx_i + \vw_i, & \quad \P(\|\vw_i\| \leq W) = 1  \\
    \widetilde{\vf}_i = \vf_i + \vw_i', & \quad \P(\|\vw_i'\| \leq W') = 1 
\end{array} $$

 The constraint~\eqref{con:lyapunov} is accordingly modified to
    $$\sum_{i=1}^{n} \; \vf_i^\top\vgamma_{i,k}+ \|\vgamma_{i,k}\|^* \left(\vphantom{\sum}M \left(\| \vy_{k,l} - \vx_{i}\| + W \right)+ W' \right)\leq s_{k,l}.$$
    For the sake of a clear presentation, we only consider noise-free measurements in the rest of this work.
\end{rem}

\subsection{Comparison with related works}~\label{sect:discussion}
We would like to wrap up this subsection by comparing the proposed learning scheme with other existing methods. In comparison with other PWA LF based methods (see e.g.~\cite{baier2012linear,marinosson2002lyapunov}), the proposed scheme shows two major differences. First, the location of the samples and the tessellation of the PWA Lyapunov candidate are decoupled in our scheme. In contrast, in existing PWA LF based methods, the data are sampled on the vertices of the tessellation (i.e. $\forall k,l$, $\exists i$ such that $\vy_{k,l} = \vx_i$), therefore, the data locations are usually structural due to the choice of the tessellation. Second, and this is a consequence of the first point, the robust Lyapunov stability conditions derived from~\eqref{eq:set_lyapunov} in the existing methods only consider the model uncertainty quantified by one data point: for instance, in~\cite{baier2012linear}, instead of our condition~\eqref{eq:robust-lyapunov}, the authors propose to prove exponential stability with the condition $\vf_i^\top\vg_k + A_k \|\vg_k\|_1 \leq - \|\vx_i\|$ for $i \in \N_n$ and some well-chosen generalization constant $A_k$ (playing the same role as our $M$), $\|\cdot\|_1$ denoting the choice $q=1$ for our norm $\|\cdot\|$. On the contrary, the scheme we propose synthetically makes use of the uncertainty quantified by each data point while maintaining a convex tractable structure.

{Another framework related to the proposed approach is the set-membership method~\cite{calliess2020lazily}. In short, the set-membership method consists in looking for a stability certificate that will be a Lyapunov function for any Lipschitz vector fields $\widehat{\vf}$ satisfying $\widehat{\vf}(\vx) \in \bfF(\vx)$ for the uncertainty set-valued map $\bfF$ from Assumption~\ref{asm:uncertainty}. In~\cite{sabug2021smgo,canale2009set}, methods are developped for the special case $d=1$ (in~\cite{canale2009set} the case $d > 1$ is deduced through assuming component-wise Lipschitz bound, i.e. $\forall j \in \N_d$, $\exists M_j > 0$ with $|\vf_j(\vx)-\vf_j(\vy)| \leq M_j \|\vx - \vy\|$), in which case $f(\vx) \in F(\vx) = [\underline{f}(\vx),\overline{f}(\vx)]$ with PWA envelope functions $\underline{f}$ and $\overline{f}$. 

Now, we will show the difference between our approach and standard set membership regarding real-valued functions. To better demonstrate the difference, we consider a specific example $f$ (Figure~\ref{fig:set_demo}), whose Lipschitz overestimate is set to $M=1$ and the data points are:
\begin{align*}
    \{\left(0,f(0)=-0.4\right),\left(0.3,f(0.3)=-0.5\right),\left(1,f(1)=-0.6\right)\}
\end{align*}
 Now, consider an LF candidate $V(x) = g \, x$ with $g=0.9$ within interval $[0,1]$. If the evaluation bounds of the set-membership method are used, the Lyapunov decreasing condition needs to be examinated in all the sub-intervals generated by the PWA bounds (plotted as two-headed arrow in Figure~\ref{fig:set_demo}). These intervals are $[0,0.1]$, $[0.1,0.2]$, $[0.2,0.3]$, $[0.3,0.6]$, $[0.6,0.7]$, $[0.7,1]$; in all generality, determining such intervals is computationally heavy. Instead, if we hope to simplify the analysis by only taking one data point into consideration as done in~\cite{baier2012linear,marinosson2002lyapunov}, none of the simple models generated by one data point can justify the Lyapunov decreasing condition, as each cone spanned by a single data point intersects the half space $gf(x) \geq 0$ (see black lines of different markers in Figure~\ref{fig:set_demo}): the use of the whole dataset is necessary. The proposed scheme synthesizes the knowledge of simple models deduced from each datapoint via convex optimization. One optimal solution to the proposed scheme is plotted as a blue line in Figure~\ref{fig:set_demo}, which only utilizes the last two points $(0.3,-0.5)$ and $(1,-0.6)$ (i.e. $\gamma_1=0$). This plot shows that the approach we propose induces a relaxation gap between the bounds it provides and the tight set-membership bounds; what we obtain in exchange for this lower accuracy is a much more tractable convex optimization scheme. It is worth noting that, in this example, if we only consider the left data and the right data point, the Lyapunov deceasing condition will fail even when the set-membership method is used. In such case, the proposed method is able to search for the additional data point locations that are relevant to the Lyapunov decreasing condition. Additionally, this process is done by polynomial time convex optimization algorithms~\cite{nesterov1994interior}. On the contrary, even though the set-membership method gives the tightest bound, checking the Lyapunov decreasing condition with these bounds is NP-hard, as it requires vertex elimination of the Voronoi cells.}

\definecolor{MyRed}{HTML}{e6550d}
\definecolor{MyBlue}{rgb}{0.20, 0.6, 0.78}
\definecolor{MyGreen}{rgb}{0.4,0.8,0.4}
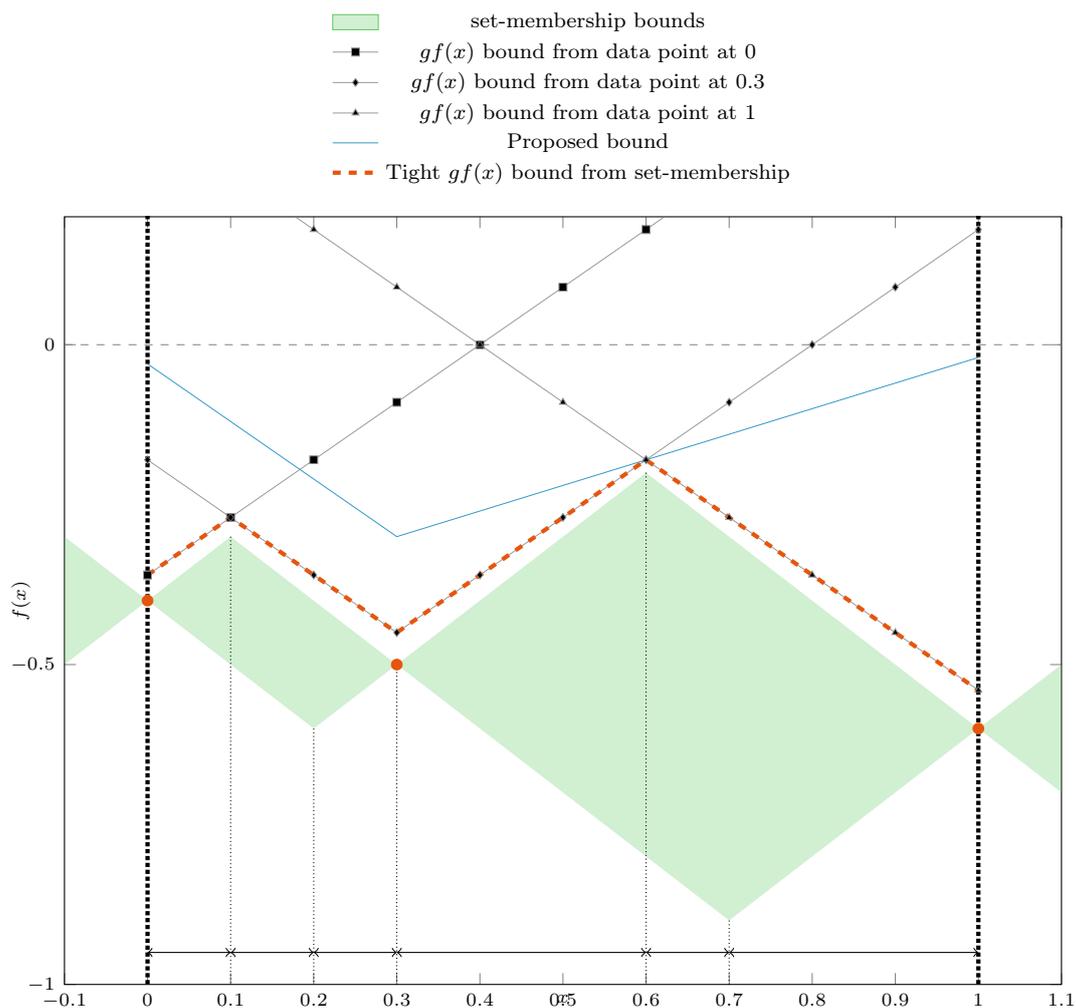
\begin{figure}[htbp!]
    \centering
    \begin{tikzpicture}
    \begin{axis}[xmin=-.1, xmax=1.1,
    ymin=-1, ymax= .2,
    enlargelimits=false,
    clip=true,
    grid=none,
    mark size=0.5pt,
    ytick distance = 0.5,
    width=1\linewidth,
    height=.8\linewidth,xlabel={$x$},ylabel = {$gf(x)$},
    legend style={
    	font=\footnotesize,
    	draw=none,
		at={(0.5,1.03)},
        anchor=south
    },
    legend columns=1,
    label style={font=\scriptsize},
    ylabel style={at={(axis description cs:0.05,0.5)}},
    xlabel style={at={(axis description cs:0.5,0.05)}},
    ticklabel style = {font=\scriptsize}]

    \pgfplotstableread{figs/compare_sm.dat}{\dat};

    \addplot+ [name path = ssub,ultra thin, MyGreen, draw opacity= 0, solid,mark=none, mark options={fill=white, scale=1.2},forget plot] table [x={x}, y ={ub}] {\dat};
    \addplot+ [name path = sslb,ultra thin, MyGreen, draw opacity= 0, solid,mark=none, mark options={fill=white, scale=1.2},forget plot] table [x={x}, y ={lb}] {\dat};
    \addplot[MyGreen,fill opacity=0.3] fill between[of=ssub and sslb]; 
    \addlegendentry{set-membership bounds}

    \addplot+ [ultra thin, gray, draw opacity= .8, solid,mark=square*, mark options={fill=black, scale=3},skip coords between index={0}{1},skip coords between index={12}{13}] table [x={x},y={one_dat1}] {\dat};
    \addplot+ [ultra thin, gray, draw opacity= .8, solid,mark=diamond*, mark options={fill=black, scale=3},skip coords between index={0}{1},skip coords between index={12}{13}] table [x={x},y={one_dat2}] {\dat};
    \addplot+ [ultra thin, gray, draw opacity= .8, solid,mark=triangle*, mark options={fill=black, scale=3},skip coords between index={0}{1},skip coords between index={12}{13}] table [x={x},y={one_dat3}] {\dat};
    \addlegendentry{$gf(x)$ bound from data point at 0}
    \addlegendentry{$gf(x)$ bound from data point at 0.3}
    \addlegendentry{$gf(x)$ bound from data point at 1}

    \addplot+ [ultra thin, MyBlue, draw opacity= 1, solid,mark=none, mark options={fill=black, scale=3},skip coords between index={0}{1},skip coords between index={12}{13}] table [x={x},y={proposed}] {\dat};
    \addlegendentry{Proposed bound}
    \addplot+ [ultra thick, dashed, MyRed, draw opacity= 1, mark=none, mark options={fill=black, scale=3},skip coords between index={0}{1},skip coords between index={12}{13}] table [x={x},y={tight_bound}] {\dat};
    \addlegendentry{Tight $gf(x)$ bound from set-membership}

    \addplot+ [ultra thin, gray, dashed,mark= none, mark options={fill=black, scale=3}] table [x={x},y={cons}] {\dat};
    
    \draw [densely dotted, ultra thick]  (axis cs: 0,1) -- (axis cs: 0, -2);
    \draw [densely dotted]  (axis cs: .1,-.3) -- (axis cs: .1, -2);
    \draw [densely dotted]  (axis cs: .2,-.6) -- (axis cs: .2, -2);
    \draw [densely dotted]  (axis cs: .3,-.5) -- (axis cs: .3, -2);
    \draw [densely dotted]  (axis cs: .6,-.2) -- (axis cs: .6, -2);
    \draw [densely dotted]  (axis cs: .7,-.9) -- (axis cs: .7, -2);
    \draw [densely dotted, ultra thick]  (axis cs: 1,1) -- (axis cs: 1, -2);

    \draw [<->] (axis cs: 0,-.95) -- (axis cs: .1, -.95);
    \draw [<->] (axis cs: .1,-.95) -- (axis cs: .2, -.95);
    \draw [<->] (axis cs: .2,-.95) -- (axis cs: .3, -.95);
    \draw [<->] (axis cs: .3,-.95) -- (axis cs: .6, -.95);
    \draw [<->] (axis cs: .6,-.95) -- (axis cs: .7, -.95);
    \draw [<->] (axis cs: .7,-.95) -- (axis cs: 1, -.95);

    \filldraw[MyRed] (axis cs: 0,-.4) circle (2pt);
    \filldraw[MyRed] (axis cs: 0.3,-.5) circle (2pt);
    \filldraw[MyRed] (axis cs: 1,-.6) circle (2pt);

    \end{axis}

    \end{tikzpicture}
    \caption{Comparison between set membership method and the proposed method. The proposed scheme is evaluated by $\Tilde{g}_1 = 0, \Tilde{g}_2 = 0.65,\Tilde{g}_3 = 0.35$.}
    \label{fig:set_demo}
\end{figure}

Recently, the set-membership approach was extended by Martin and Allgöwer in~\cite{martin2023poly}. In their contribution, the authors extend Lipschitz bounds to derivative bounds (see Example in Section~\ref{sect:stage}), that are used to deduce a functional set-membership parameterization, namely they look for certificates (for dissipativity rather than attractivity) that are consistent with any vector field $\widehat{\vf} \in \cF_\bfD^{\vomega}$, for some approximation point $\vomega \in \R^d$. The key idea in their methodology is that $\cF_\bfD^{\vomega}$ can be represented in terms of polynomial sums of squares, that are themselves parameterized by linear matrix inequalities (LMI). This very general framework competes with other set membership methods in terms of accuracy, and successfully rids itself of the corresponding NP-hardness by approximating positivity constraints with LMIs through the celebrated Positivstellensatz. However, one of the involved LMIs has size $2\binom{d+k}{d} + d$ (see $\Xi_i$ in LMI (8) of~\cite{martin2023poly}), i.e. combinatorial in the number of state variables and the maximal order of the derivative bounds, hence preventing the method from scaling to high dimensions. In contrast, the approach we propose relies on SOCP and LP formulations, that are known to scale much better than SDP.

To summarize this discussion, the approach we propose can be seen as a trade-off between computationally efficient methods based on PWA LFs such as~\cite{baier2012linear,marinosson2002lyapunov} on one hand, and more accurate set-membership based methods as~\cite{calliess2020lazily,sabug2021smgo,canale2009set,martin2023poly} on the other hand. As a result, it achieves better accuracy and sample-efficiency than the former (see the example illustrated on Figure~\ref{fig:set_demo}) while being theoretically more scalable than the latter.

\section{Algorithm Development}\label{sect:num_details}
After the introduction of the Lyapunov learning problem~\ref{pb:lasalle}, we will discuss its learnability in Section~\ref{sect:valid}. The original learning problem will be recast to an equivalent but numerically more efficient form in Section~\ref{sect:comp_imprv}. In the end, the main algorithms are summarized in Section~\ref{sect:algo}.

\subsection{Learnability}\label{sect:valid}

First, due to the introduction of the slack variables $\{s_{k,l}\}$, the learning problem~\ref{pb:lasalle} is always feasible. It is natural to ask the following core question in the limiting case: 

\begin{center}
    \textit{Given an RoA prior estimate $\bfA$, what condition on dataset $\bfD$ should hold to enable learning the PWA Lyapunov candidate on $\bfX\setminus\bfA$?}
\end{center}

Obviously, it is impossible to answer this question with a sufficient condition, regarding the arbitrariness of the unknown dynamic system $\vf$. However, if we only assume a bound on the Lipschitz constant of $\vf$ as in Assumption~\ref{asm:lipschitz}, we can still give an initial check of the learnability of problem~\ref{pb:lasalle}, based on the remark we already made on the fact that $\bfF(\vx)$ should not contain $\vO$ for $\vx \in \bfX\setminus \bfA$. In order to discuss this necessary condition, we first define 
\begin{align*}
    \rho_i := \frac{\lVert \vf_i\rVert}{M}\;.
\end{align*}
We state the necessary condition as follows:
\begin{lem}\label{lem:covering}
With tessellation~\eqref{eq:tessellation} and under Assumptions~\ref{asm:regularity},~\ref{asm:data},~\ref{asm:prior},~\ref{asm:poly} and~\ref{asm:lipschitz}, if solutions to problem~\ref{pb:lasalle} define an LF as in Proposition~\ref{prop:lasalle}, then it holds
\begin{equation} \label{eq:covering}
    \bfX\setminus\bfA\subset \bigcup_{i=1}^{n}\bfB(\vx_i,\rho_i).
\end{equation}
\end{lem}
\begin{proof}
As we already proved, under our assumptions, the existence of an optimal solution to problem~\ref{pb:lasalle} satisfying condition~\eqref{eq:lasalle} implies that the Lyapunov condition~\eqref{eq:pwa-lyapunov} is satisfied on $\bfX\setminus\bfA$, i.e.
$$ \forall k \in \N_m, \; \vx \in \bfY_k, \; \vf \in \bfF(\vx), \quad \vf^\top \vg_k < 0. $$
In particular, for all $k \in \N_m$, $\vx \in \bfY_k$, it holds $\vO \notin \bfF(\vx)$, otherwise one would get the contradiction $0 = \vO^\top\vg_k < 0$. Now, from Assumption~\ref{asm:lipschitz} we know that 
$$ \bfF(\vx) = \bigcap_{i=1}^n\bfB(\vf_i, M\|\vx - \vx_i\|), $$
so that $\vO \notin \bfF(\vx)$ means existence of an $i \in \N_n$ such that $\vO \notin \bfB(\vf_i , M\|\vx-\vx_i\|)$, i.e.
$$ M\rho_i = \|\vf_i\| = \|\vf_i - \vO\| > M\|\vx - \vx_i\|,$$
and hence $\vx \in \bfB(\vx_i, \rho_i)$. This yields the announced result:
$$ \bfX \setminus \bfA \subset \bigcup_{k=1}^m \bfY_k \subset \bigcup_{i=1}^n \bfB(\vx_i, \rho_i).$$
\end{proof}

\begin{rem}
    Again, notice that the necessary condition~\eqref{eq:covering} for learnability implies the Borel-Lebesgue property, i.e. that $\overline{\bfX \setminus \bfA}$ is compact, and hence that $\bfX$ is bounded, which further supports our assumption in that direction.
\end{rem}

Lemma~\ref{lem:covering} shows the connection between learning a PWA LF and the set covering problem, which was proved to be equivalent to a non-convex semi-infinite problem~\cite{krieg2019modeling}, and thus one should not try to check the condition in Lemma~\ref{lem:covering} numerically. Recall a key idea behind the problem~\ref{pb:lasalle}: the global analysis on $\bfX\setminus\bfA$ is reduced to the analysis on the vertices. This inspires us to relax the continuous set covering problem to the covering problem of the vertices and we state this condition in the following Corrolary:
\begin{cor}\label{cor:covering}
With tessellation~\eqref{eq:tessellation} and under Assumptions~\ref{asm:regularity},~\ref{asm:data},~\ref{asm:prior},~\ref{asm:poly} and~\ref{asm:lipschitz}, if solutions to problem~\ref{pb:lasalle} define an LF as in Proposition~\ref{prop:lasalle}, then it holds
 \begin{equation}\label{eqn:tess_valid}
    \{\vy_{k,l}\}_{\substack{1\leq k\leq m \\ 1\leq l \leq \nu_k}} \subset \bigcup_{i=1}^{n}\bfB(\vx_i,\rho_i).
\end{equation}
\end{cor}

This necessary condition~\eqref{eqn:tess_valid} can be checked in polynomial time. If this test fails, it means that there exist a $\vy_{k,l} \notin \cup_{i=1}^{n}\bfB(\vx_i,\rho_i)$, and therefore, the dataset $\bfD$ is not informative enough to learn a PWA LF by only assuming the Lipschitz condition~\eqref{eq:lipschitz}. Accordingly, the learning process will be aborted. Intuitively, the vertices $\vy_{k,l} \notin \cup_{i=1}^{n}\bfB(\vx_i,\rho_i)$ should suggest the location where additional samples are required. We leave the investigation about this aspect for future work.

\subsection{Computationally efficient recasting}\label{sect:comp_imprv}
In this part, we will discuss how we recast the original problem~\ref{pb:lasalle} to an equivalent problem that can be handled numerically more efficiently.

\noindent\textbf{Data Refinement}\\
One main computational bottleneck for the original problem~\ref{pb:lasalle} comes from the number of decision variables. Without loss of generality, we consider an affine piece $\bfY_k$ and look for indices $i$ that could be removed from $\N_n$ without jeopardizing the tessellation validation condition~\eqref{eqn:tess_valid}. In other words, we look for $i \in \N_n$ such that, for all $l \in \N_{\nu_k}$, $\vy_{k,l} \notin \bfB(\vx_i,\rho_i)$. For such index $i$, and for any $\vgamma \in \R^d$, it holds $\|\vx_i - \vy_{k,l}\| > \rho_i$ and hence, using Hölder's inequality:
\begin{align*}
    \vf_i^\top\vgamma + M \|\vgamma\|^*\|\vx_i-\vy_{k,l}\| & > \vf_i^\top\vgamma + M \rho_i \|\vgamma\|^* \\
    & = \vf_i^\top\vgamma + \|\vgamma\|^*\|\vf_i\| \\
    & \geq \vf_i^\top\vgamma - \vgamma^\top\vf_i = 0.
\end{align*}

Hence, giving a nonzero value to the $\vgamma_{i,k}$ corresponding to this choice of $k$ and $i$ cannot help enforcing the strict negative Lyapunov condition~\eqref{con:lyapunov}, so that it is optimal to fix such $\vgamma_{i,k}$ to $\vO$. Such choice is equivalent to removing the index $i$ from the sum in the Lyapunov condition~\eqref{con:lyapunov}
$$ \sum_{i=1}^n \vf_i^\top\vgamma_{i,k} + M\|\vgamma_{i,k}\|^*\|\vx_i-\vy_{k,l}\| \leq s_{k,l}. $$
Hence, we can define the set of data indices relevant to $\bfY_k$: 
\begin{subequations} \label{eq:indices}
\begin{equation} \label{eq:relevant_indices}
    \I_k := \left\{i \in \N_n \; \middle| \; \bfB(\vx_i,\rho_i) \cap \{\vy_{k,l}\}_{l=1}^{\nu_k} \neq \varnothing\right\}
\end{equation}
and replace the dataset $\bfD = \{(\vx_i,\vf_i)\}_{i=1}^n$ with $\bfD_k := \{(\vx_i,\vf_i)\}_{i\in\I_k} \subset \bfD$ when working on the polytopic cell $\bfY_k$, recasting condition~\eqref{con:lyapunov} as
\begin{equation} \label{eq:relevant-lyapunov}
    \sum_{i\in\I_k} \vf_i^\top\vgamma_{i,k} + M\|\vgamma_{i,k}\|^*\|\vx_i-\vy_{k,l}\| \leq s_{k,l}.
\end{equation}
\end{subequations}

The condition defining the new dataset $\bfD_k$ essentially states that each datapoint should at least bring useful information to enforce strict negativity on constraint~\eqref{con:lyapunov} at one vertex. This technique can significantly reduce the computational cost. To see this, we consider a homogeneous tessellation within a unit hypercube centered at $\mathbf{0}$ within which data points scatter uniformly. We further assume that the cells of the tessellation are hypercubes with edge width $\eta$. Then, considering a datapoint $\vx_i$, the ball $\bfB(\vx_i,\rho_i)$ intersects at most $O((\nicefrac{\rho_i}{\eta})^d) \ll m$ hypercubic cells. 
Notice that in that specific case, it holds 
$$\rho_i = \frac{\|\vf_i\|}{M} = \frac{\|\vf_i-\vO\|}{M} \leq \frac{M\|\vx_i-\vO\|}{M} = \|\vx_i\| \leq 0.5,$$ 
so that the number of decision variables are reduced to roughly $O\left(\frac{1}{2^{d}}\right)$ of the problem defined by the whole data set. In the numerical example we consider in Section~\ref{sect:result}, we observe on average an 81\% reduction in the number of decision variables, which makes the problem tractable on a Laptop without memory overflow.\\

\noindent\textbf{Explicit SOCP formulations}\\
When dealing with the euclidean norm (i.e. $q=2$), in the numerical implementation, it is critical to convert the inequality constraint~\eqref{eq:relevant-lyapunov} into a set of second order cone constraints~\cite{boyd2004convex}:
\begin{subequations}
\begin{align}
    s_{k,l}\geq& \sum_{i \in \I_k} \; \vf_i^\top \vgamma_{i,k} + z_{i,k,l} \label{eq:socp-slack} \\
    z_{i,k,l} &\geq M \lVert \vgamma_{i,k}\rVert \lVert \vy_{k,l}-\vx_{i}\rVert, \label{eq:socp-lyapunov}
\end{align}
\end{subequations}
where $|\I_k|$ auxiliary decision scalar variables $\{z_{i,k,l}\}$ are introduced per vertex $\vy_{k,l}$. More precisely, the complexity for solving an SOCP scales like $O(N+K+L)^{3.5}$ where $N$ is the number of variables, $K$ the cumulated dimension of the cone constraints, and $L$ the number of linear constraints, while solving an SDP scales like $O(N+K^2+L)^{3.5}$, where $K$ denotes the size of the involved LMI. In our case, $N$, $K$ and $L$ all scale linearly with the sample size $n$, resp. state dimension $d$, resp. number of tessellation vertices $\nu = \nu_1+\ldots+\nu_m$, so that the overall complexity would be $O(P^7)$ ($P \in \{n,d,\nu\}$ being the varying parameter) without the SOCP reformulation, versus $O(P^{3.5})$ with the SOCP recasting.

\noindent\textbf{The Recast Problem}\\
After the above reformulations, one ends up numerically solving the following problem:
\begin{pb} \label{pb:socp}
Find $\{\{s_{k,l}\}_{l=1}^{\nu_k}\}_{k=1}^m$, $\{\{z_{i,k,l}\}_{1 \leq l \leq \nu_k}^{i\in\I_k}\}_{k=1}^{m}$, $\{b_k\}_{k=0}^m \subset \R$, $\{\{\vgamma_{i,k}\}_{i\in\I_k}\}_{k=1}^{m} \subset \R^d$, $\vg_0 \in \R^d$ solution to
\begin{subequations} \label{opt:socp}
\begin{align}
& s_\epsilon^\star := \min \sum\limits_{k=1}^m \sum\limits_{l=1}^{\nu_k} s_{k,l} \notag \\
& \text{if } \vy_{k,l} \in \bfY_0, \quad \sum_{i \in \I_k} \vgamma_{i,k}^\top \vy_{k,l} = \vg_0^\top\vy_{k,l} + b_0 - b_k \tag{\ref*{eqn:pce-aff-cont}} \\
& \text{if } \vy_{k,l} \in \bfY_{k'}, \quad \sum_{i \in \I_k} (\vgamma_{i,k}-\vgamma_{i,k'})^\top \vy_{k,l} = b_{k'}-b_k \tag{\ref*{eq:decomposition}} \\
& \forall k \in \N_m,\ i\in \I_k,\ l \in \N_{\nu_k}, \quad s_{k,l} \geq -\epsilon \tag{\ref*{con:slack-lasalle}} \\
& \qquad \quad s_{k,l}\geq \sum_{i \in \I_k} \; \vf_i^\top \vgamma_{i,k} + z_{i,k,l} \tag{\ref*{eq:socp-slack}} \\
& \text{and } \quad z_{i,k,l} \geq M \lVert \vgamma_{i,k}\rVert^* \lVert \vy_{k,l}-\vx_{i}\rVert. \tag{\ref*{eq:socp-lyapunov}}
\end{align}
\end{subequations}
\end{pb}

\begin{rem}
    \textbf{Learning without knowing $M$: }It is possible to consider $M$ as a decision variable, and  determine the largest Lipschitz constant for which an LF can be found given a particular set of data. More specifically, constraint~\eqref{eq:socp-lyapunov} is recast to a positive semi-definite constraint by 
    \begin{align}\label{robust-lyapunov_sdp}
        \begin{bmatrix}
        MI_{d}&\vgamma_{i,k}\\\vgamma_{i,k}^\top & z_{i,k,l}
        \end{bmatrix}\in\mathcal{S}^+\;,
    \end{align}
    where $\mathcal{S}^+$ denotes the set of positive semi-definite matrices. However, once this formulation is used, the resulting optimization problem becomes an SDP and the data refinement technique proposed at the beginning of this Section~\ref{sect:comp_imprv} cannot be applied. Even though the particular sparsity structure in~\eqref{robust-lyapunov_sdp} can be exploited to improve the computational efficiency, the computational cost of this optimization still drastically increase in comparison with the SOCP problem~\ref{pb:socp}.
\end{rem}

\subsection{Algorithms} \label{sect:algo}
A first learning scheme is summarized in Algorithm~\ref{alg:learn_Lya}. Even though we use a standard, non-tailored tessellation methodology in algorithm~\ref{alg:learn_Lya}, generating a good tessellation is vital but non-trivial. Existing works mostly focus on the link between a convex liftable tessellation and the power diagram (see e.g.~\cite{rybnikov2000polyhedral,aurenhammer1987power}). However, as the LF studied in this paper is not necessarily convex, hence we leave the study of this topic in the future research and we use the standard Delaunay triangulation in this work~\cite{delaunay1934sphere}.
\begin{algorithm}
\caption{$\enspace$}\label{alg:learn_Lya}
\hspace*{\algorithmicindent} \textbf{Input}: RoA prior $\bfA$, negativity tolerance $-\epsilon$, Lipschitz overestimate $M$ \newline
\hspace*{\algorithmicindent} \textbf{Output}: LF $V(\vx)$, RoA inner estimate $\widehat{\bfR} = \LV_a$

\begin{algorithmic}
    \STATE Refine a tessellation $\{\bfY_k\}_{k=1}^{m}$ until it satisfies~\eqref{eqn:tess_valid}
    \IF{tessellation is valid}
        \STATE Solve optimization problem~\ref{pb:socp}
        \IF{Optimal solution satisfies $s^\star_{\epsilon}=-\epsilon\sum_{k=1}^{m} \nu_k$}
            \STATE Find the largest $a$ such that $\LV_a \subset \bfX$
            \STATE Return RoA estimate $\widehat{\bfR} = \LV_a$
        \ENDIF
    \ELSE
        \STATE Return cannot learn $V(\vx)$.
    \ENDIF 
\end{algorithmic}
\end{algorithm}

\noindent\textbf{Proving that $\bfX$ is the full RoA}
Up to this point, we only assumed knowledge of a subset $\bfA$ of the RoA and focused on finding an RoA estimate based on the dataset $\bfD$ and the Lipschitz bound $M$, $\bfX$ being the admissible state set in which the considered system should evolve. However, in some cases, one could define $\bfX$ as the a priori candidate for the RoA estimate. In that case, one wants to certify that the whole $\bfX$ is equal to the RoA $\bfR$, by finding a PWA LF $V(\vx)$ and level value $a \in \R$ such that $\bfX = \LV_a$. This is done by adding the following level conditions to problem~\ref{pb:socp} (recalling that $\bfX$ is \textit{open}, i.e. $\bfX \cap \partial \bfX = \varnothing$):
\begin{subequations} \label{eq:level}
\begin{align} 
    \vy_{k,l} \in \partial \bfX \implies \quad & \sum_{i\in\I_k} \vgamma_{i,k}^\top\vy_{k,l} + b_k = a \label{eq:boundary} \\
    \vy_{k,l} \in \bfX \implies \quad & \sum_{i\in\I_k} \vgamma_{i,k}^\top\vy_{k,l} + b_k \leq a - \epsilon. \label{eq:interior}
\end{align}
\end{subequations}

Based on the solution to the resulting problem, it is possible to determine an RoA estimate using the following Corollary:
\begin{cor}
    Under Assumption~\ref{asm:poly} and tessellation given by~\eqref{eq:tessellation}, for any $\vf$ satisfying Assumptions~\ref{asm:regularity},~\ref{asm:data},~\ref{asm:prior}, and~\ref{asm:lipschitz}, if the solution to problem~\ref{pb:socp} further satisfies satisfies~\eqref{eq:lasalle} i.e. $s^\star_{\epsilon}=-\epsilon\sum_{k=1}^{m} \nu_k$, as well as level condition~\eqref{eq:level} for some $a \in \R$, then $\bfX = \bfA$ i.e.
    $$ \forall\vx_0 \in \bfX, \qquad \lim_{t\to\infty} \|\vx(t|\vx_0)\| = 0. $$
\end{cor}
\begin{proof}
    Trivially follows from Proposition~\ref{prop:lasalle} applied to $\LV_a \stackrel{\eqref{eq:level}}{=} \bfX$.
\end{proof}
With this corollary, we are able to certify whether $\bfX$ is the full RoA, by applying the scheme summarized in Algorithm~\ref{alg:learn_Lya_noX}.

\begin{algorithm}
\caption{$\enspace$}\label{alg:learn_Lya_noX}
\hspace*{\algorithmicindent} \textbf{Input}: RoA prior $\bfA$, negativity tolerance $-\epsilon$, Lipschitz overestimate $M$, set level $a$ \newline
\hspace*{\algorithmicindent} \textbf{Output}: LF $V(\vx)$

\begin{algorithmic}
    \STATE Refine a tessellation $\{\bfY_k\}_{k=1}^{m}$ until it satisfies~\eqref{eqn:tess_valid}
    \IF{tessellation is valid}
        \STATE Solve optimization problem~\ref{pb:socp} with constraint~\eqref{eq:level}
        \IF{Optimal value solution satisfies $s^\star_{\epsilon}=-\epsilon\sum_{k=1}^{m} \nu_k$}
            \STATE Return
        \ENDIF
    \ELSE
        \STATE Return cannot learn $V(\vx)$.
    \ENDIF 
\end{algorithmic}
\end{algorithm}

It is possible to sequentially decompose the learning of the RoA: in addition to increasing the size of the RoA prior $\bfA$, we propose to sequentially increase the size of the RoA candidate $\bfX$, i.e. we come up with a sequence $(\bfX_r)_{1\leq r\leq r_{\max}}$ of RoA candidates, such that

$$\bfA \subset \bfX_1 \subset \bfX_2 \ldots \subset \bfX_{r_{\max}} = \bfX. $$

Now, at iteration $r$, we prove that $\bfX_r$ is positively invariant and attracted to $\vO$, after which we set $\bfA = \bfX_r$ and work on $\bfX_{r+1}$; to do so, it is only required to learn the LF on the set $\bfX_r \setminus \bfA$, as summarised in Algorithm~\ref{alg:learn_seq_noX}. At the end of the iteration process, and if the algorithm did not fail, one gets an algorithmic proof that $\bfX_{r_{\max}} = \bfX = \bfR$. If the algorithm failed at iteration $r$, then one still gets an algorithmic proof that $\bfX_{r-1} \subset \bfR$ (with convention $\bfX_0 = \bfA$ at the first iteration). It is noteworthy that, if one needs to recover the whole LF on $\bfX\setminus\bfA$, then it is necessary to impose the continuity condition on the boundary of $\partial\bfX_r$ between the $r$-th iteration and the $(r+1)$-th iteration, as the tessellation of $\bfX \setminus \bfA$ is the union of tessellations of the $\bfX_{r+1} \setminus \bfX_{r}$, so that the continuity condition~\eqref{eqn:pce-aff-cont} applies on the boundary $\partial\bfX_r$.

\begin{algorithm} 
\caption{$\enspace$}\label{alg:learn_seq_noX}
\hspace*{\algorithmicindent} \textbf{Input}: RoA prior $\bfA$, RoA candidates $(\bfX_r)_{1\leq r \leq r_{\max}}$ \newline
\hspace*{\algorithmicindent} \textbf{Output}: LF $V(\vx)$, RoA inner estimate $\widehat{\bfR} = \LV_a$

\begin{algorithmic}
\STATE $r=1$
\WHILE{$\{r \leq r_{\max}\}$}
    \STATE Run Algorithm.~\ref{alg:learn_Lya_noX} with $\bfX = \bfX_r$
    \IF{Algorithm~\ref{alg:learn_Lya_noX} failed}
        \STATE Return cannot learn $V(\vx)$.
    \ENDIF 
    \STATE $\bfA \gets \bfX_r$
    \STATE $r \gets r+1$
\ENDWHILE
\end{algorithmic}
\end{algorithm}

\section{Numerical Results}\label{sect:result}
In this part, the proposed learning schemes are evaluated in two different examples. In particular,  we will make use of Algorithm~\ref{alg:learn_Lya},~\ref{alg:learn_Lya_noX} and ~\ref{alg:learn_seq_noX} in the first example and use Algorithm~\ref{alg:learn_Lya} in the second one. All the following results are implemented on a laptop with Intel i7-11800H and 32G memory, and the \texttt{Mosek} numerical solver is used to get optimal solutions of the SOCP problems.

\subsection{Non Polynomial dynamic system}
We consider a two-dimensional nonlinear dynamic system:
\begingroup\makeatletter\def\f@size{9.5}\check@mathfonts
\begin{align*}
    \dot{x}_1(t) =& -0.9\sin(x_1(t))\cos(x_2(t))+0.2x_1(t)x_2(t)+0.25x_2(t)^2\\
    \dot{x}_2(t) =& -1\sin(x_2(t))(|x_1(t)+0.2|)+0.5\frac{x_1(t)x_2(t)}{\cos(x_2(t))-0.3x_1(t)}\\
\end{align*}
\endgroup
\begin{figure}[h!]
\centering
\begin{subfigure}{0.49\textwidth}
    \centering
    \includegraphics[width=0.9\textwidth]{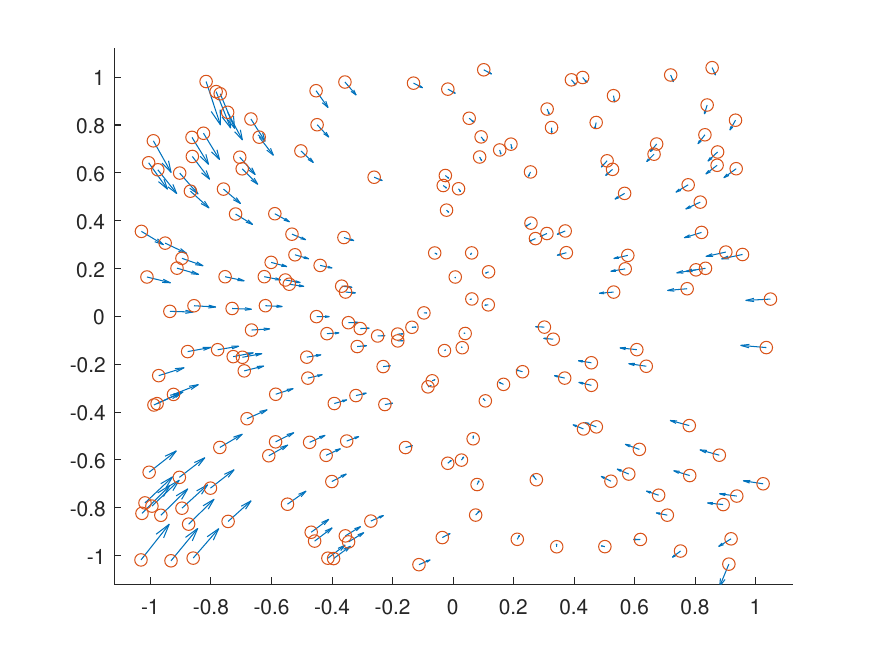}
    \caption{The data used for the learning scheme.}
    \label{fig:vis_sample}
\end{subfigure}
\begin{subfigure}{0.49\textwidth}
    \centering
    \includegraphics[width=0.9\textwidth]{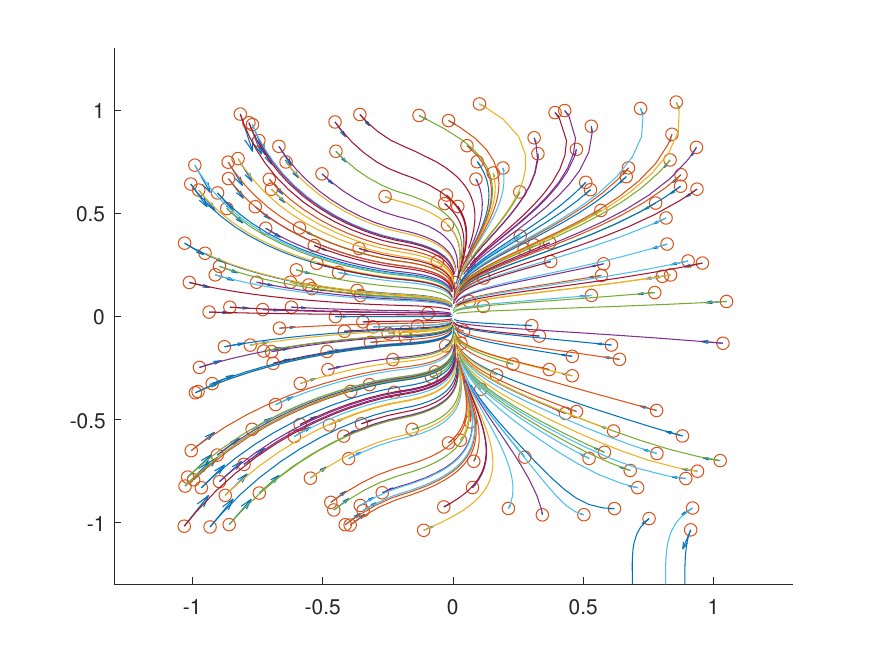}
    \caption{Underlying dynamics (not used for learning).}
    \label{fig:traj}
\end{subfigure}
\caption{Representations of the considered nonpolynomial dynamic system.}
\end{figure}

We assume that we know an RoA prior $\bfA=[-0.1,0.1]^2$. A dataset with only 200 samples within $[-1,1]^2\subset\mathbb{R}^2$ is used to learn the underlying LF: the positions $\{\vx_i\}$ and speeds $\{\vf_i\}$ of these samples are plotted in Figure~\ref{fig:vis_sample}, from which we can observe that this dataset is relatively sparse in $[-1,1]^2$. Judging by the speed sample, the dynamic system seems stable within the box $[-0.4,0.4]^2$, while stability within the region $[-1,1]^2\setminus [-0.4,0.4]^2$ is unclear because of the speed samples in the lower right corner in Figure~\ref{fig:vis_sample}. Hence, we ran sequential space partition scheme (Algorithm~\ref{alg:learn_seq_noX}). In particular, we first use Algorithm~\ref{alg:learn_Lya} in the region $[-0.4,0.4]^2$ with $\bfA=[-0.1,0.1]^2$. After we justify that $[-0.4,0.4]^2$ is a positively invariant subset of the RoA, then we further apply Algorithm~\ref{alg:learn_Lya_noX} to $[-1,1]^2$ with $\bfA=[-0.4,0.4]^2$. In both sub-problems, the negativity tolerances $\epsilon$ are set to $10^{-3}$ and the tessellation are both randomly generated by Delaunay triangulation~\cite{delaunay1934sphere}.

\begin{figure}[h!]
\centering
\begin{subfigure}{0.49\textwidth}
    \centering
    \includegraphics[width=0.9\textwidth]{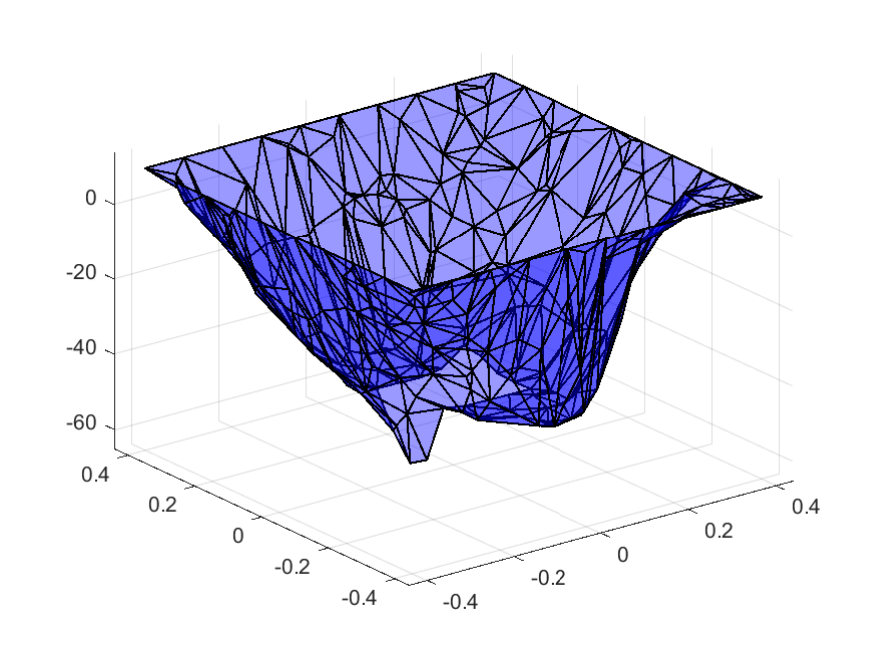}
    \caption{$\cX_s = [-0.1,0.1]^2 \subset \cX = [-0.4,0.4]^2$.}
    \label{fig:lya1}
\end{subfigure}
\begin{subfigure}{0.49\textwidth}
    \centering
    \includegraphics[width=0.9\textwidth]{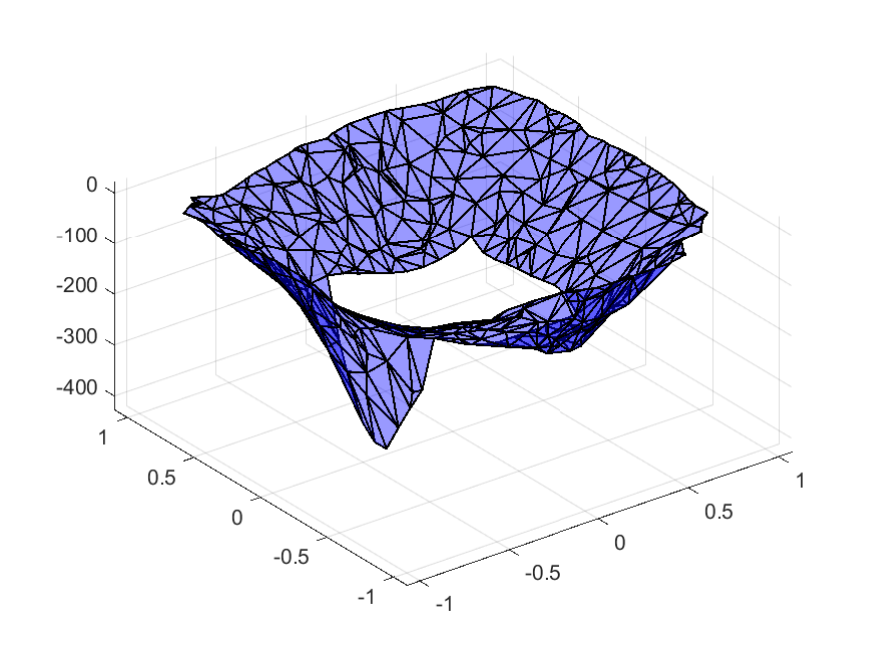}
    \caption{Learning a RoA estimate within $[-1,1]^2$.}
    \label{fig:lya2}
\end{subfigure}
	\caption{Visualization of the learnt LFs.}
\end{figure}
\begin{figure}[h!]
    \centering
    \includegraphics[width=.5\textwidth]{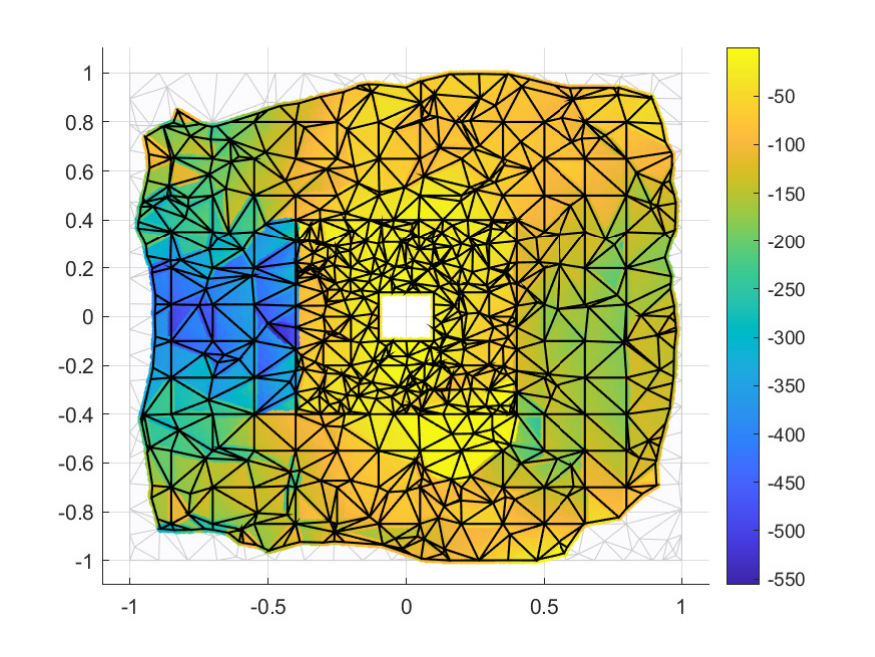}
    \caption{Evaluation of $\vf(x)^\top\vpartial_{\rm C}V(\vx)$ on the learnt RoA, the gray triangulation in the background is the tessellation used to solve Problem~\ref{pb:socp}, while the coloured region in the front is our RoA estimate.}
    \label{fig:eval}
\end{figure}

The learnt Lyapunov function in $[-0.4,0.4]^2$ is shown in Figure~\ref{fig:lya1}, while the RoA we finally end up with is shown in Figure~\ref{fig:eval}. Moreover, the LF learnt from Algorithm~\ref{alg:learn_Lya_noX} in $\bfA=[-1,1]^2\setminus [-0.4,0.4]^2$ is shown in Figure~\ref{fig:lya2}. Figure~\ref{fig:eval} also shows the evaluation of $\vf(x)^\top\vpartial_{\rm C}V(\vx)$ with respect to the underlying dynamic system, whose maximal evaluation is $-1.525\times10^{-2} < 0$, as expected from an LF. In accordance with our guess, the learnt RoA in Figure~\ref{fig:eval} cuts off the lower right corner, because this region does not seem to be stable. To see that, we simulate the underlying dynamic system by setting the initial states to points in our dataset. The simulated trajectories are plotted in Figure~\ref{fig:traj}; please note that these trajectories are not used in the learning scheme at all. 

\subsection{Reverse Time Van Del Pol Oscillator}

In this part, we consider the reverse time Van Del Pol oscillator:

\begin{align*}
    \dot{x}_1(t) &= -2x_2(t)\\
    \dot{x}_2(t) &= -0.8x_1(t)-10(x_1(t)^2-0.21)x_2(t)\;.
\end{align*}

\begin{figure}[h!]
\centering
\begin{subfigure}{0.49\textwidth}
    \centering
    \includegraphics[width=0.9\textwidth]{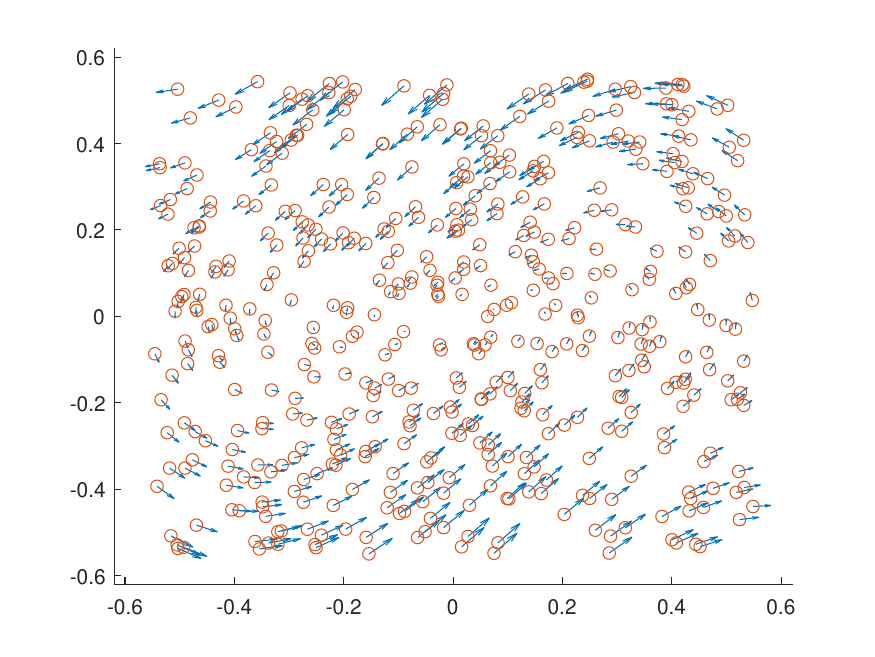}
    \caption{The data used for the learning scheme.}
    \label{fig:vis_sample_vdp}
\end{subfigure}
\begin{subfigure}{0.49\textwidth}
    \centering
    \includegraphics[width=0.9\textwidth]{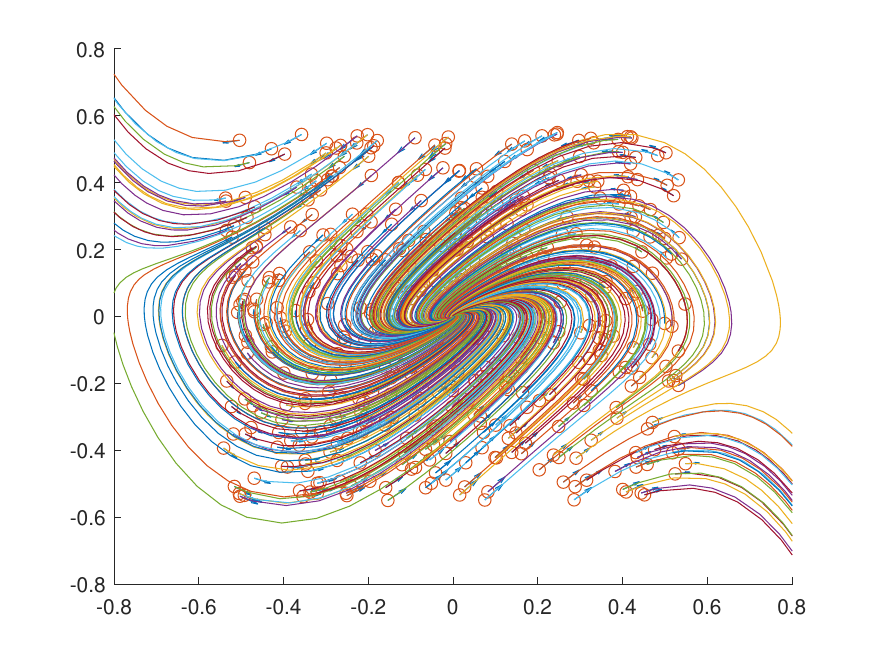}
    \caption{Simulated trajectories (not used for learning).}
    \label{fig:traj_vdp}
\end{subfigure}
\caption{The reverse time Van Del Pol oscillator.}
\end{figure}

\begin{figure}[h!]
\centering
\begin{subfigure}{0.49\textwidth}
    \centering
    \includegraphics[width=0.9\textwidth]{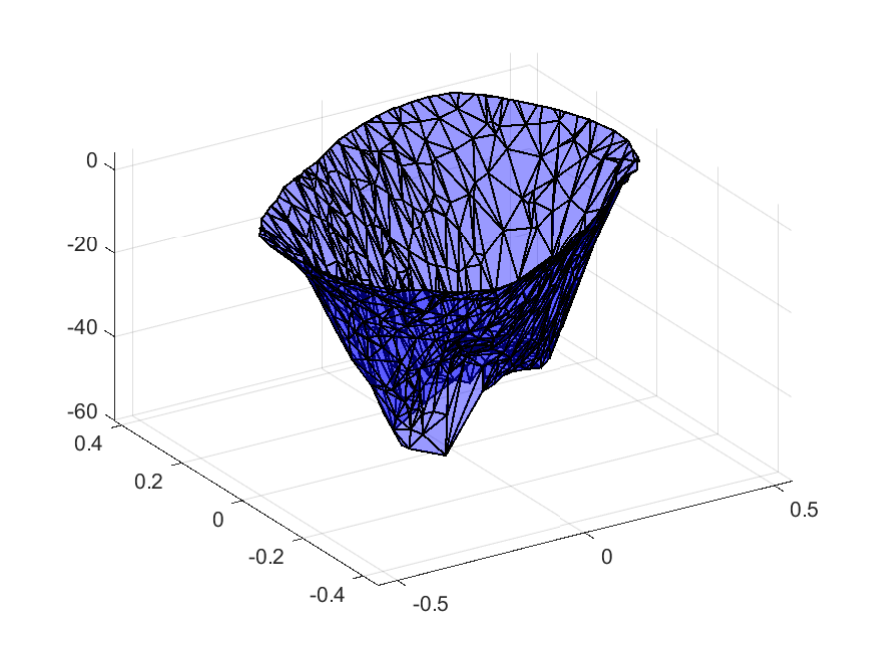}
    \caption{Visualization of the LF.}
    \label{fig:lya_vdp}
\end{subfigure}
\begin{subfigure}{0.49\textwidth}
    \centering
    \includegraphics[width=0.9\textwidth]{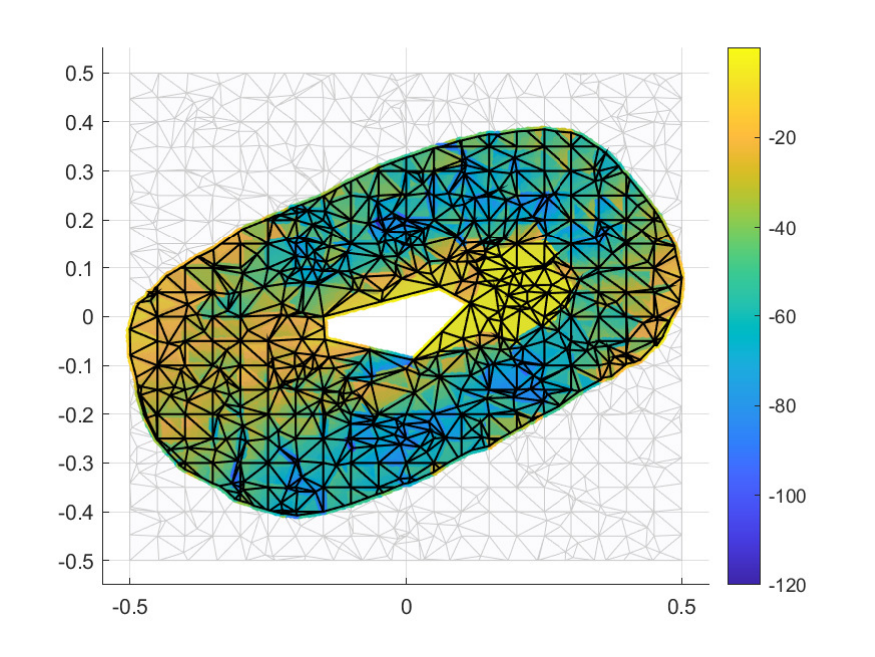}
    \caption{Evaluation of $\vf(\vx)^\top\vpartial_{\rm C}V(\vx)$ on the learnt RoA.}
    \label{fig:eval_vdp}
\end{subfigure}
\caption{The Lyapunov function learnt from a Van der Pol dataset.}
\flushleft
On Figure~\ref{fig:eval_vdp}, the gray triangulation in the background is the tessellation used to solve Problem~\ref{pb:socp}, while the coloured region in the front is the RoA estimate. The polytopic hole in the middle represents the RoA prior $\bfA$.
\end{figure}

We know an a-priori polytopic RoA $\bfA$, which is plotted in the center of Figure~\ref{fig:eval_vdp}. A dataset with only 400 samples within $[-0.5,0.5]^2\subset\mathbb{R}^2$ is used to learn the underlying LF: the positions $\{\vx_i\}$ and speeds $\{\vf_i\}$ of these samples are plotted in Figure~\ref{fig:vis_sample_vdp}. Similar to what we did in the last example, we simulate these data forward in Figure~\ref{fig:traj_vdp}, while these trajectories are not used in the learning scheme. We can observe that both the lower right corner and the upper left corner in Figure~\ref{fig:traj_vdp} correspond to regions of unstable states. Even with only the access to the data in Figure~\ref{fig:vis_sample_vdp}, we can not give a clear idea about which region is safe, hence we apply Algorithm~\ref{alg:learn_Lya_noX} to  $[-0.5,0.5]^2$. In particular, the negativity tolerance $\epsilon$ is set to $10^{-3}$ and the tessellation is randomly generated by Delaunay triangulation~\cite{delaunay1934sphere}. The learnt LF and the corresponding evaluation of $\vf(\vx)^\top\vpartial_{\rm C}V(\vx)$ on the learnt RoA are respectively plotted in Figure~\ref{fig:lya_vdp} and Figure~\ref{fig:eval_vdp}. In particular, the maximal evaluation of $\vf(\vx)^\top\vpartial_{\rm C}V(\vx)$ on the learnt RoA is $-1.947\times 10^{-2} < 0$.

\section{Conclusion}\label{sect:conclusion}

The paper proves a variant of a stability theorem with non-smooth Lyapunov functions (LF) and then implements an algorithm for data-based region of attraction (RoA) estimation of an unknown dynamic system. Through this process, a theorem for piecewise affine (PWA) LF computation was proven, robustly distributing the global dataset information and deriving a convex optimization program for computing such Lyapunov functions.

The originality of the method is that it only requires a fixed dataset to compute an estimate of the RoA, from which it allows the user to deduce global information from local data and knowledge of the Lipschitz constant of the unknown dynamic system. In more detail, the technique (1) decouples the PWA representation from the data point locations, and (2) decomposes the corresponding uncertainty set containing the unknown dynamics as an intersection of sets, each one depending on a single data point, in order to derive a tractable criterion accounting for the global dataset. Hence, it can be used to study systems whose dynamics cannot be easily sampled at will, through a relatively simple optimization problem that can be handled with interior point methods.

\bibliographystyle{ieeetr}
\bibliography{ref}

\appendix

\section{Proof of Theorem~\ref{thm:rb_stb_Lip}} \label{app:proof}

\textit{Theorem~\ref*{thm:rb_stb_Lip}}: 
Let Assumptions~\ref{asm:regularity},~\ref{asm:data},~\ref{asm:prior},~\ref{asm:poly} and~\ref{asm:lipschitz} hold, and consider the function $V$ defined by~\eqref{eq:tessellation},~\eqref{eq:cpwalf}. For $k \in \N_m$, let $\{\vy_{k,l}\}_{l=1}^{\nu_k}$ be the set of vertices of $\bfY_k$, so that 
\begin{equation} \tag{\ref*{eq:vertices}}
\bfY_k = \co\{\vy_{k,l}\}_{l=1}^{\nu_k}.
\end{equation}
Suppose that $\forall k \in \N_m$,~there is a set $\{\vgamma_{i,k}\}_{i=1}^n \subset \R^d$ with
\begin{equation} \tag{\ref*{eq:deco-Lip}}
\sum_{i=1}^n \vgamma_{i,k} = \vg_k
\end{equation}
and for all $l \in \N_{\nu_k}$
\begin{equation}\tag{\ref*{eq:robust-lyapunov}}
\sum\limits_{i=1}^{n} \vf_i^\top\vgamma_{i,k}+ M \|\vgamma_{i,k}\|^*\| \vx_i - \vy_{k,l}\| <0.
\end{equation}

Then, $\forall a \in \R$ such that $\LV_a \subset \bfX$, it holds $\LV_a \subset \bfR$, i.e.
\begin{equation} \tag{\ref*{eq:stability}}
\forall \vx_0 \in \LV_a, \qquad \lim_{t\to\infty}\|\vx(t|\vx_0)\| = 0.
\end{equation}
\begin{proof}
    By Assumption~\ref{asm:lipschitz}, for all $\vx \in \bfX$, $i\in\N_n$, it holds $\|\vf(\vx) - \vf_i\| \leq M \|\vx - \vx_i\|$, i.e. $\vf(\vx) \in \bfB(\vf_i, M\|\vx-\vx_i\|) =: \bfF_i(\vx)$. In particular, 
$$\vf(\vx) \in \bfF(\vx) = \bigcap_{i=1}^n \bfF_i(\vx),$$
so that Assumptions~\ref{asm:uncertainty} and~\ref{asm:localize} hold. Hence, according to Proposition~\ref{prop:stability_pwa} and Theorem~\ref{thm:robust}, with $\bOm = \{\vx_i\}_{i=1}^n$, the proof of our statement boils down to verifying~\eqref{eq:distribution}, i.e.
$$ \sum_{i=1}^n \max\limits_{\vf \in \bfF_i(\bfY_k)} \vf^\top \vgamma_{i,k} < 0. $$
Let $i \in \N_n$, $k \in \N_m$, $\vf \in \bfF_i(\bfY_k)$. 
By construction, there exists an $\vx \in \bfY_k$ such that $\vf \in \bfF_i(\vx)$, i.e. $\|\vf-\vf_i\|\leq M \|\vx-\vx_i\|$, and the H\"{o}lder inequality $\forall \vphi,\vpsi\in\R^d, |\vphi^\top\vpsi| \leq \|\vphi\|^*\|\vpsi\|$, yields
\begin{align*}
\vf^\top\vgamma_{i,k} & = \vf_i^\top\vgamma_{i,k} + \vgamma_{i,k}^\top(\vf - \vf_i) \\
& \leq \vf_i^\top\vgamma_{i,k} + \|\vgamma_{i,k}\|^*\|\vf - \vf_i\| \\
& \leq \vf_i^\top\gamma_{i,k} + M \|\vgamma_{i,k}\|^*\|\vx-\vx_i\|.
\end{align*}
Taking the maximum over $\vf \in \bfF_i(\bfY_k)$, one gets
$$ \max_{\vf \in \bfF_i(\bfY_k)} \vf^\top\vgamma_{i,k} \leq \vf_i^\top\vgamma_{i,k} + M \|\vgamma_{i,k}\|^* \left(\max_{\vx \in \bfY_k} \|\vx_i-\vx\| \right). $$
It remains to show that
$$\max_{\vx \in \bfY_k} \|\vx_i-\vx\| = \max_{1\leq l \leq \nu_k} \|\vx_i - \vy_{k,l}\|, $$
which follows from the convexity of the objective function $\vx \longmapsto \|\vx_i-\vx\|$ together with the convexity of $\bfY_k$: the maximum of the optimization problem on the left hand side is attained on an extreme point of $\bfY_k$, i.e. a $\vy_{k,l}$. Finally, summing over $i$ yields
\begin{align*}
& \sum_{i=1}^n \max_{\vf\in\bfF_i(\bfY_k)} \vf^\top \vgamma_{i,k} \leq \max_{1\leq l \leq \nu_k} \sum_{i=1}^n \vf_i^\top\vgamma_{i,k} \\
& \qquad \quad + M \|\vgamma_{i,k}\|^* \| \vx_i - \vy_{k,l} \| \stackrel{\eqref{eq:robust-lyapunov}}{<} 0.
\end{align*}
\end{proof}

\begin{rem} If instead of Assumption~\ref{asm:lipschitz}, we suppose that $\|\vf\|_{\cH^d} \leq M$ for some RKHS $\cH$ with kernel $\kappa$, then the reasoning given in the Example of Section~\ref{sect:stage} yields an uncertainty set 
$$ \bfF_i(\vx) = \frac{\kappa(\vx,\vx_i)}{\kappa(\vx_i,\vx_i)}\vf_i + \underset{P_i(\vx)}{\underbrace{\sqrt{\kappa(\vx,\vx) - \frac{\kappa(\vx,\vx_i)^2}{\kappa(\vx_i,\vx_i)}}}} \left[\pm \vmu_i \right] $$
for some fixed vector $\vmu_i \in \R^d$ depending on $M$ and the representer $\widehat{\vf}$ of $\vf$ learnt from the single data $\{\vx_i,\vf_i\}$ (but not on $\vx$), so that in the previous proof one has to replace $M \|\vx-\vx_i\|$ with $ |\nicefrac{\kappa(\vx,\vx_i)}{\kappa(\vx_i,\vx_i)} - 1| \cdot \|\vf_i\| + 2 P_i(\vx) \|\vmu_i\|$ (using the triangle inequality to bound $\|\vf-\vf_i\|$ for $\vf \in \bfF_i(\vx)$ as above). However, the problem here is that this new function of $\vx$ is not convex: its maximum is not necessarily attained on a vertex $\vy_{k,l}$ of $\bfY_k$, and we would need other arguments (out of the scope of this article) to obtain a tractable constraint.
\end{rem}

\section{Assessing conservatism}

One can see from the proof of Theorem~\ref{thm:robust} that replacing~\eqref{eq:set_lyapunov} with~\eqref{eq:robust} induces some conservatism. This subsection is devoted to proving that, in the special case of Lipschitz bounds, this conservatism ultimately vanishes when the data points and tessellation vertices form an appropriate covering of $\bfX \setminus \bfA$.

\begin{prop}
    Let $V : \R^d \longrightarrow \R$ be a continuously differentiable LF with $V(\vO)=0$, $V > 0$ and $\dot{V}_{\vf} < 0$ on $\bfX \setminus \{\vO\}$ (which is proved to exist under Assumption~\ref{asm:regularity}) and suppose that there is an $a \in \R$ with $\bfA \subset \LV_a \subset \bfX$.

    Under Assumptions~\ref{asm:data},~\ref{asm:prior},~\ref{asm:poly} and~\ref{asm:lipschitz} and up to adding samples to our dataset $\bfD$, $V$ can be approximated by a continuous PWA function $\hat{V}$ whose parameters in the sense of~\eqref{eq:cpwalf} are solution to Problem~\ref{pb:lasalle} and consistent with condition~\eqref{eq:lasalle}, so that $\hat{V}$ defines a generalized LF as per Proposition~\ref{prop:lasalle}.
\end{prop}
\begin{proof}
    Let $\delta > 0$. Using the universal approximation property of PWA functions~\cite[Chapter 7.4]{royden1988real}, there exists a tessellation~\eqref{eq:tessellation} and a PWA function $\hat{V}$ as in~\eqref{eq:cpwalf} such that for all $k \in \{0,\ldots,m\}$, $\vx \in \bfY_k$, it holds
    \begin{subequations} \label{eq:proof}
    \begin{equation} \label{eq:approx}
        |V(\vx) - \vg_k^\top\vx - b_k| < \delta, \qquad |\dot{V}_{\vf}(\vx) - \vf(\vx)^\top\vg_k| < \delta.
    \end{equation}
    We choose $\delta > 0$ small enough such that 
    \begin{equation} \label{eq:beta}
        \beta := -\max_{\vx \in \overline{\bfX \setminus \bfA}} \dot{V}_{\vf}(\vx) - \delta > 0.
    \end{equation}
    Let $\eta > 0$. Up to adding samples to our dataset $\bfD$ and further partitioning the tessellation cells $\bfY_k$, we assume that
        \begin{align}
            & \bfX \setminus \bfA \subset \bigcup_{i=1}^n \bfB(\vx_i, \nicefrac{\eta}{M}) \label{eq:ref_data} \\
            & \forall k \in \N_m, \quad \max_{\vx,\vy \in \bfY_k} \|\vx-\vy\| < \frac{\eta}{M} \label{eq:ref_tess}
        \end{align}
    Let $k \in \N_m$ and $\vx \in \bfY_k$. By~\eqref{eq:ref_data}, there exists an $i_k \in \N_n$ such that $\|\vx_{i_k} - \vx\| < \nicefrac{\eta}{M}$. As a result, it holds
    \begin{align}
        \vf_{i_k}^\top\vg_k & \leq \vf(\vx)^\top\vg_k + \|\vg_k\|^*\|\vf_{i_k} - \vf(\vx) \| \notag \\
        & \stackrel{\eqref{eq:approx}}{\leq} \dot{V}_{\vf}(\vx) + \delta + \|\vg_k\|^*\|\vf_{i_k} - \vf(\vx) \| \notag \\
        & \stackrel{\eqref{eq:beta}}{\leq} \|\vg_k\|^*\|\vf_{i_k} - \vf(\vx) \| - \beta \notag \\
        & \stackrel{\eqref{eq:lipschitz}}{\leq} M\|\vg_k\|^*\|\vx_{i_k} - \vx\| - \beta \notag \\
        & \stackrel{\eqref{eq:ref_data}}{\leq} \eta\|\vg_k\|^* - \beta. \label{eq:eval}
    \end{align}
    We use this to define
    $$ \vgamma_{i,k} = \begin{cases}
        \vg_k \quad \text{if } i=i_k \\
        \vO \quad \text{else.}
    \end{cases} $$
    Then, it holds for all $l \in \N_{\nu_k}$
    \begin{align}
        \sum_{i=1}^n \vf_i^\top\vgamma_{i,k} + M\|\vgamma_{i,k}\|^*\|\vy_{k,l}-\vx_i\| \hspace*{-12em} & \notag \\
        & = \vf_{i_k}^\top\vg_k + M\|\vg_k\|^*\|\vy_{k,l}-\vx_{i_k}\| \notag \\
        & \stackrel{\eqref{eq:eval}}{\leq} \|\vg_k\|^*\left(\eta + M\|\vy_{k,l}-\vx_{i_k}\|\right) - \beta \notag \\
        & \stackrel{(*)}{\leq} 3\eta\|\vg_k\|^* - \beta,\label{eq:proof_bound}
    \end{align}
    \end{subequations}
    where $(*)$ is obtained through the simple decomposition
    \begin{align*}
        \|\vy_{k,l}-\vx_{i_k}\| & \leq \|\vy_{k,l} - \vx\| + \|\vx-\vx_{i_k}\| \\
        & \stackrel{\eqref{eq:ref_data}}{\leq} \|\vy_{k,l} - \vx\| + \nicefrac{\eta}{M} \stackrel{\eqref{eq:ref_tess}}{\leq} 2\nicefrac{\eta}{M}
    \end{align*}
    In particular, setting
    $$ \eta = \frac{\epsilon + \beta}{3\max_{k\in\N_m}\|\vg_k\|^*} $$
    yields, for all $k\in\N_m$, $l \in \N_{\nu_k}$,
    $$ \sum_{i=1}^n \vf_i^\top\vgamma_{i,k} + M\|\vgamma_{i,k}\|^*\|\vy_{k,l}-\vx_i\| \stackrel{\eqref{eq:proof_bound}}{\leq} -\epsilon $$
    so that our choice for $\{b_k\}_{k=0}^m$, $\{\vgamma_{i,k}\}_{1\leq i\leq n}^{1\leq k \leq m}$, $\vg_0$ and $\{\{s_{k,l}\}_{l=1}^{\nu_k}\}_{k=1}^m = \{-\epsilon\}$ yields a feasible solution to Problem~\ref{pb:lasalle}.
\end{proof}

\end{document}